\documentclass[11pt,a4paper]{article}
\usepackage[T2A]{fontenc}
\usepackage[cp1251]{inputenc}
\usepackage[english]{babel}
\usepackage{epsfig,amssymb,epic,eepic,color, graphicx}
\usepackage{fancyhdr}
\usepackage{amsthm}
\usepackage{amsmath}
\usepackage{indentfirst}
\usepackage{float}
\usepackage{enumitem}   

\usepackage{colortbl}

\usepackage{pb-diagram}

\usepackage{amsbib}

\usepackage{authblk}

\usepackage{setspace}
\usepackage{float}
\usepackage{fancyhdr}
\usepackage{geometry}
\newtheorem{theorem}{Theorem}
\newtheorem{definition}{Definition}
\newtheorem{corollary}{Corollary}
\newtheorem{remark}{Remark}
\newtheorem{lemma}{Lemma}[section]
\newtheorem{proposition}{Proposition}[section]
\newenvironment{demo}{{\bf Proof: }}{\hfill $\square$ \medskip}

\newenvironment{demo2}{{\bf Proof of Theorem 2: }}{\hfill $\square$ \medskip}
\newenvironment{demo3}{{\bf Proof of Theorem 3: }}{\hfill $\square$ \medskip}
\newenvironment{demo4}{{\bf Proof of Theorem 4: }}{\hfill $\square$ \medskip}

\newenvironment{democor1}{{\bf Proof of Corollary 1: }}{\hfill $\square$ \medskip}

\renewenvironment{abstract}{\textbf{Abstract.} }{\medskip}
\begin{document}

\begin{center}
\Large{High order homoclinic tangencies of corank 2}
\end{center}

\begin{center}
\large{Dmitrii Mints}
\end{center}

\begin{center}
\large{Department of Mathematics, Imperial College London, London, UK}
\end{center}

\begin{abstract}
We prove that in the space of $C^r$ maps $(r=2,\ldots,\infty,\omega)$ of a smooth manifold of dimension at least 4 there exist open regions where maps with infinitely many corank-2 homoclinic tangencies of all orders are dense. The result is applied to show the existence of maps with universal two-dimensional dynamics, i.e. maps whose iterations approximate the dynamics of every map of a two-dimensional disk with an arbitrarily good accuracy. We show that maps with universal two-dimensional dynamics are $C^r$-generic in the regions under consideration.
\end{abstract}

\section{Introduction}\label{section1}

The problem of classification for dynamical systems is solvable for systems with a uniformly hyperbolic chain-recurrent set (the set of such systems coincides with the set of the Axiom A systems \cite{Smale} which satisfy the no-cycle condition \cite{stability1},\cite{stability2}; this includes Anosov systems, Morse-Smale systems, horseshoe maps, etc.). Hyperbolic systems are structurally stable in the sense that if two such systems are $C^1$ close, then they are topologically equivalent in a neighborhood of the chain-recurrent set. Hyperbolic systems comprise an open set in the space of smooth systems, and the structural stability entails that they are described by a discrete set of invariants. 

The major fact in the theory of dynamical systems is that the complement to the set of structurally stable systems has a non-empty interior -- the regions of structural instability, and the dynamics for systems from these regions are extremely diverse and much more complicated than in the hyperbolic case. The simplest mechanism of destroying the structural stability is a \textit{homoclinic tangency}, i.e. a non-transverse intersection of the stable and unstable manifolds of a hyperbolic periodic orbit. Although any given tangency between two manifolds is a fragile object, the presence of homoclinic tangencies in a system turns out to be a persistent phenomenon. By the Newhouse Theorem \cite{Newhouse2},\cite{Newhouse3}, there exist $C^2$-open regions (the \textit{Newhouse domain}) in the space of dynamical systems where systems with homoclinic tangencies are dense. Importantly, Newhouse regions are found in the space of parameters of many popular systems of applied interest, including the H{\'e}non map, the Chua circuit, the Lorenz and R{\"o}ssler models, etc. 

According to \cite{Turaevintro1},\cite{Turaevintro2}, the main characterisitc of systems from the Newhouse domain is the ultimate richness of the dynamics. In line with this thesis, we show in the present paper that there exist $C^2$-open subregions of the Newhouse domain, where systems with extremely degenerate homoclinic tangencies are dense. The result implies the existence of a new type of universal dynamics which are generic for these regions.

\subsection{Highly degenerate tangencies}\label{subsection1.1}

The tangency between two manifolds is itself a degenerate property (a generic intersection of two manifolds is transverse). We characterize the degree of the degeneracy by the \textit{corank} and the \textit{order} of the tangency. Let manifolds $W_1$ and $W_2$ form a tangency at a point $M$. Let $\mathcal T_M W$ denotes the tangent space to a manifold $W$ at $M$. 

\begin{figure}[h]
\begin{center}
\begin{minipage}[h]{0.47\linewidth}
\includegraphics[width=1\linewidth]{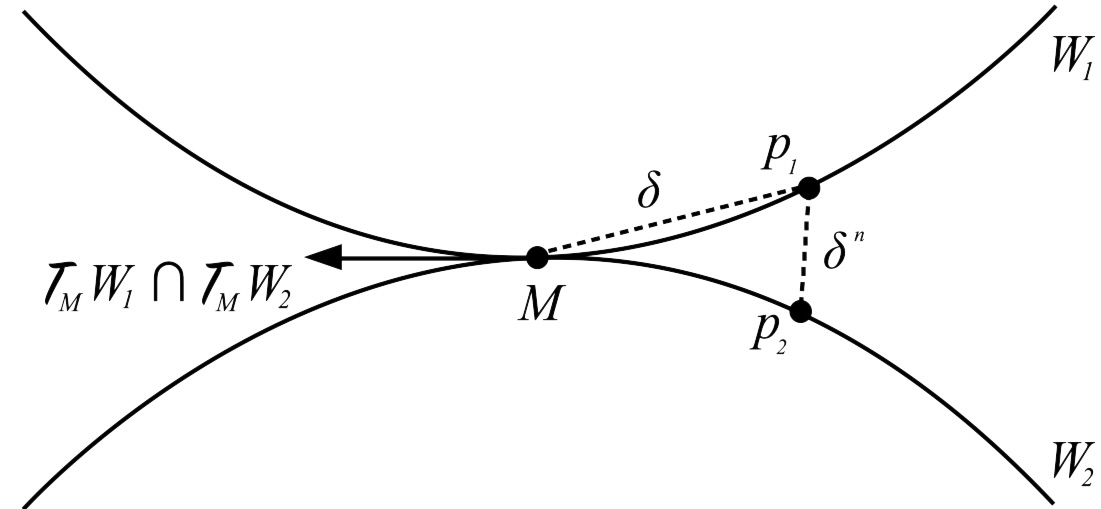}
\caption{Curves $W_1$ and $W_2$ have a corank-1 tangency of order $n$ at the point $M$.} 
\label{Figure1} 
\end{minipage}
\hfill 
\begin{minipage}[h]{0.47\linewidth}
\includegraphics[width=1\linewidth]{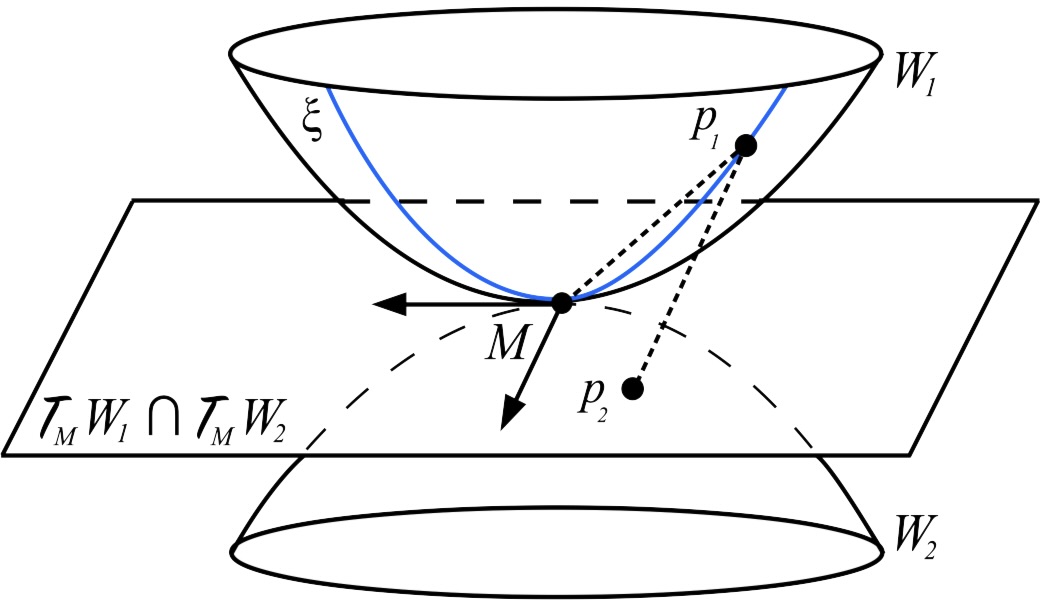}
\caption{Manifolds $W_1$ and $W_2$ have a corank-2 tangency at the point $M$.}
\label{Figure2}
\end{minipage}
\end{center}
\end{figure}

We say that \textit{the tangency is of corank} $\mathcal C$ if  
\begin{equation*}
\begin{aligned} 
&\dim \mathcal T_M W_1+\dim \mathcal T_M W_2 -\dim\left(\mathcal T_M W_1\oplus\mathcal T_M W_2\right)=\dim\left(\mathcal T_M W_1\cap \mathcal T_M W_2\right)=\mathcal C.
\end{aligned}
\end{equation*}
If $W_1$ is one-dimensional (a curve), then we say that \textit{the tangency} with $W_2$ at the point $M$ \textit{is of order $n$} if for any point $p_1\in W_1$, which is sufficiently close to $M$, there exists a point $p_2\in W_2$ such that
\begin{equation*}
\begin{aligned} 
&dist(p_1,p_2)=O\left(dist(p_1,M)^n\right).
\end{aligned}
\end{equation*}
If $\dim W_1>1$, let $\Xi$ be the set which consists of smooth curves $\xi\subset W_1$ such that $M\in\xi$ and $\mathcal T_M \xi \subseteq \mathcal T_M W_1\cap \mathcal T_M W_2$. Then we denote the order of the tangency between a curve $\xi\in\Xi$ and $W_2$ as $\text{ord}_{\xi}$, and say that \textit{the tangency} between $W_1$ and $W_2$ at the point $M$ \textit{is of order $n$} if
\begin{equation*}
\begin{aligned} 
&\min\limits_{\xi\in\Xi} \text{ord}_{\xi}=n
\end{aligned}
\end{equation*}
(see Fig. \ref{Figure1} and \ref{Figure2}). Thus, the order of tangency measures how flat the tangency is: the higher the order, the flatter the tangency (see a coordinate definition in Section \ref{subsection3.2}). 

A generic tangency is of corank 1, that is, the manifolds $W_1$ and $W_2$ are tangent along one direction only. Systems with corank-1 homoclinic tangencies have been studied extensively over the past decades, and currently a well-developed theory exists (see references in the books \cite{Palis2},\cite{Bonatti1},\cite{homoclinictangencies}). The central result (the Newhouse Theorem) is that  in the space of smooth dynamical systems, near any system with a homoclinic tangency there exist $C^2$-open regions where systems with corank-1 homoclinic tangencies are dense (see \cite{Newhouse2},\cite{Newhouse3},\cite{Gonchenko1},\cite{Palis1},\cite{Romero1},\cite{Gourmelon}). 

One of the most important properties of the Newhouse domain is the density of maps with corank-1 homoclinic tangencies of \textit{arbitrarily high orders} \cite{arbitraryorder}, \cite{arbitraryorderconservative}. More precisely, maps with \textit{infinitely many homoclinic tangencies of all orders} are dense in the Newhouse domain in the space of two-dimensional $C^r$ maps ($r\ge 2$; this includes $C^{\infty}$, and the space of real-analytic maps). This result has many applications, including universal dynamics for area-preserving maps \cite{arbitraryorderconservative} and Beltrami fields \cite{Berger1}. It is also essential for constructing maps with highly degenerate periodic orbits \cite{arbitraryorder}, \cite{arbitraryorderconservative}, maps with superexponential growth of the number of periodic points\footnote{Other mechanisms for obtaining superexponential growth of the number of periodic points are discussed in \cite{periodicpoints1},\cite{periodicpoints2},\cite{periodicpoints3}.} \cite{Kaloshin1} and "unsmoothable" diffeomorphisms (which cannot be topologically conjugate to any diffeomorphism of a higher regularity) \cite{highsmoothness}. In the multidimensional case, the density (in the Newhouse domain) of maps having infinitely many corank-1 homoclinic tangencies of all orders is proven in \cite{Mints1}.

Systems with tangencies of high corank ($\mathcal C>1$) turn out to be a complicated object to describe, and their study began only a few years ago \cite{Barrientos1},\cite{Barrientos2},\cite{Barrientos3},\cite{Buzzi1},\cite{Catalan1},\cite{Asaoka}. In the present paper, we concentrate on the phenomenon of high order tangencies of corank 2. Our main result is 

\begin{theorem}\label{theorem1}
In the space of $C^r$ maps $(r=2,...,\infty,\omega)$ of each $k$-dimensional manifold with $k\ge 4$ there exist open regions where maps with infinitely many orbits of corank-2 homoclinic tangencies of every order form a dense subset.
\end{theorem}

Further, we discuss the main idea of the proof and explain the construction of open regions from Theorem \ref{theorem1}. 

\subsubsection{Main idea} 
 
Let $f$ be a $C^r$ diffeomorphism $(r=2,\ldots,\infty,\omega)$ of a smooth or analytic $k$-dimensional manifold $\mathcal M^k$. Suppose that $f$ has a periodic orbit $L$ of period $h$, that is, $L=\{O,f(O),...,f^{h-1}(O)\}$ with $f^h(O)=O$ and $f^i(O)\neq O$ for all $0<i<h$. The eigenvalues of the Jacobian matrix for the map $f^h$ calculated at the point $O$ are called \textit{multipliers} of the periodic orbit $L$. The orbit $L$ is called \textit{hyperbolic} if none of its multipliers lies on the unit circle. Any hyperbolic periodic orbit is structurally stable in the sense that if $L_f$ is such an orbit of the map $f$, then any map $g$, which is $C^r$-close\footnote{Here and throughout the rest of the article we measure the distance between two maps on some compact region $K\subset\mathcal M^k$. We say that two $C^r$-smooth maps are $\delta$-close if a $C^r$ distance between them on $K$ does not exceed $\delta$. If $r=\infty$, we define a $C^{\infty}$ distance as $\rho_{\infty}(f_1,f_2)=\sum\limits^{\infty}_{r=0}\frac{1}{(r+1)^2}\cdot\frac{\rho_r(f_1,f_2)}{1+\rho_r(f_1,f_2)}$, where $\rho_r$ is a $C^r$ distance. If $r=\omega$ (the real-analytic case), we fix some small complex neighborhoods $Q$ of $K$ and say that two $C^{\omega}$ maps are $\delta$-close if they differ on not more than $\delta$ at every point of $Q$.} to $f$, has a hyperbolic periodic orbit $L_g$ which is a continuation of $L_f$.

We will assume that the hyperbolic periodic orbit $L$ has multipliers on both sides of the unit circle. Let $\lambda_1,...,\lambda_{k_s}$, $\gamma_1,\ldots,\gamma_{k_u}$ be multipliers of $L$ ordered so that $|\gamma_{k_u}|\ge...\ge|\gamma_1|>1>|\lambda_1|\ge\ldots\ge|\lambda_{k_s}|$. The multipliers inside the unit circle are said to be stable and those outside the unit circle are said to be unstable. Denote $\lambda=|\lambda_1|$, $\gamma=|\gamma_1|$. Those multipliers which are equal in absolute value to $\lambda$ or $\gamma$ are called \textit{leading}, and the rest are called \textit{non-leading}. For the orbit $L$ there exist a smooth $k_1$-dimensional \textit{stable manifold} $W^s(L)$ and a smooth $k_2$-dimensional \textit{unstable manifold} $W^u(L)$ defined as
\begin{equation}\label{stableandunstablemanifolds}
\begin{aligned} 
&W^s(L)=\{x\in\mathcal M^k: dist(f^m(x),L)\to 0 \;\; \text{as} \;\; m\to+\infty\},\\
&W^u(L)=\{x\in\mathcal M^k: dist(f^{-m}(x),L)\to 0 \;\; \text{as} \;\; m\to+\infty\}.
\end{aligned}
\end{equation}

We will call the orbit $L$ a \textit{bi-focus}, if both its leading stable and leading unstable multipliers are complex conjugate. So, $\lambda_{1,2}=\lambda e^{\pm i\varphi}, \; \gamma_{1,2}=\gamma e^{\pm i\psi}$, where $\varphi,\psi\not=0,\pi$. This condition implies that the ambient manifold dimension $k$ is at least 4. 

The main idea in the proof of Theorem \ref{theorem1} is an algorithm which allows one to obtain a high order homoclinic tangency of corank 2 by adding a $C^r$-small perturbation to a system with a given (sufficiently large) number of quadratic homoclinic tangencies of corank 2 formed by the stable and unstable manifolds of a bi-focus periodic orbit. It is an independent result that can be applied to solving other problems (see the discussion in Section \ref{subsection1.3}), therefore we formulate it as a separate

\begin{theorem}\label{theorem3}
Let $f$ be a $C^r$ map $(r=2,\ldots,\infty,\omega)$ with a bi-focus periodic orbit $L_f$ whose stable and unstable manifolds contain $2^{\frac{(N-1)(N+4)}{2}}$ orbits of a corank-2 homoclinic tangency. Then arbitrarily $C^r$-close to $f$ there exists a map $g$ with a bi-focus periodic orbit $L_g$ whose stable and unstable manifolds contain an orbit of a corank-2 homoclinic tangency of order $N$.
\end{theorem}
 
\subsubsection{$ABR^*$-domain}\label{subsubsection1.1.2}
 
Let us recall that a compact, topologically transitive, uniformly hyperbolic and locally maximal\footnote{A set $\Lambda$ is called locally maximal if there exists its compact neighborhood $V$ such that $\Lambda=\bigcap\limits_{i\in\mathbb Z} f^{-i}(V)$. We will call such a neighborhood defining.} set $\Lambda$ of a smooth map $f$ is called a \textit{basic set}. If a basic set is just a single hyperbolic periodic orbit, then it is said to be trivial; otherwise it is called non-trivial. Throughout this paper, all basic sets are zero-dimensional. 

Let $\Lambda$ be a basic set of the map $f$ and let $p$ be a point in $\Lambda$. Define \textit{stable} $W^s(p)$ and \textit{unstable} $W^u(p)$ \textit{manifolds} of a point $p\in\Lambda$ as
\begin{equation*}
\begin{aligned} 
&W^s(p)=\{x\in\mathcal M^k: dist(f^m(x),f^m(p))\to 0 \;\; \text{as} \;\; m\to+\infty\},\\
&W^u(p)=\{x\in\mathcal M^k: dist(f^{-m}(x),f^{-m}(p))\to 0 \;\; \text{as} \;\; m\to+\infty\}.
\end{aligned}
\end{equation*}
The topological dimension of the stable and unstable manifolds is the same for all points of the basic set $\Lambda$, therefore one can define the \textit{type} $(k_s,k_u)$ of $\Lambda$ where $k_s$ and $k_u$ are the dimensions of the stable and unstable manifolds, respectively. The manifolds $W^s(p)$ and $W^u(p)$ are injective $C^r$-immersions of $\mathbb R^{k_s}$ and $\mathbb R^{k_u}$, respectively,  and they depend continuously on the point $p$ and on the map $f$.

In the same way as this is done for a hyperbolic periodic orbit (which is a trivial basic set; see \eqref{stableandunstablemanifolds}), one can introduce \textit{stable} $W^s(\Lambda)$ and \textit{unstable} $W^u(\Lambda)$ \textit{manifolds} of an arbitrary basic set $\Lambda$ as
\begin{equation*}
\begin{aligned} 
&W^s(\Lambda)=\{x\in\mathcal M^k: dist(f^m(x),\Lambda)\to 0 \;\; \text{as} \;\; m\to+\infty\},\\
&W^u(\Lambda)=\{x\in\mathcal M^k: dist(f^{-m}(x),\Lambda)\to 0 \;\; \text{as} \;\; m\to+\infty\}.
\end{aligned}
\end{equation*}
Let us note that 
\begin{equation*}
\begin{aligned} 
&W^s(\Lambda)=\bigcup\limits_{p\in\Lambda} W^s(p), \;\;\;\; \text{and} \;\;\;\; W^u(\Lambda)=\bigcup\limits_{p\in\Lambda} W^u(p).
\end{aligned}
\end{equation*}

According to \cite{Smale}, \cite{Bowen}, the sets $W^s(L)\cap\Lambda$ and $W^u(L)\cap\Lambda$ are dense in the basic set $\Lambda$ for any periodic orbit $L\in\Lambda$. Moreover, stable and unstable manifolds of $L$ accumulate to stable and unstable manifolds of points from $\Lambda$ in the following sense: for any point $p\in\Lambda$ and any compact ball $K_p$ in $W^s(p)$ or $W^u(p)$ there exists a compact ball $K_L$ in $W^s(L)$, respectively in $W^u(L)$, such that $K_p$ and $K_L$ are $C^r$-close. 

Let $\Lambda_f=\bigcap\limits_{i\in\mathbb Z} f^{-i}(V)$, where $V$ is a defining neighborhood of $\Lambda_f$, be a non-trivial basic set of type $(k_s,k_u)$ of the map $f$. It is well-known that for every map $g$ which is $C^r$-close to $f$ the set $\Lambda_g=\bigcap\limits_{i\in\mathbb Z} g^{-i}(V)$ (the continuation of $\Lambda_f$) is also a non-trivial basic set of type $(k_s,k_u)$. We say that $\Lambda_f$ \textit{exhibits a $C^r$-robust tangency of corank $\mathcal C$} if there exists a compact subset $\mathcal L\subset W^s(\Lambda_f)\cup W^u(\Lambda_f)$ with the following property: given any compact neighborhood $U$ of $\mathcal L$ there exists a $C^r$ neighborhood $\mathcal U$ of $f$ such that for any map $g\in\mathcal U$ there are points $p_1, p_2$ in the continuation $\Lambda_g$ for which $W^s(p_1)$ and $W^u(p_2)$ have an orbit of corank-$\mathcal C$ tangency completely contained in $U$. 

In recent papers by Barrientos and Raibekas \cite{Barrientos1} and Asaoka \cite{Asaoka}, the following result is proved.

\begin{proposition}\label{propositionhighcorank}\cite{Barrientos1},\cite{Asaoka}
Every manifold of dimension $k\ge 4$ admits a map with a non-trivial basic set exhibiting a $C^2$-robust tangency of corank $\mathcal C$, which can be chosen to be any integer $0<\mathcal C\le[\frac{k}{2}]$.
\end{proposition}

Let us emphasize that the bound on $\mathcal C$ in Proposition \ref{propositionhighcorank} is sharp since for a basic set of type $(k_s,k_u)$ the stable and unstable manifolds have dimension $k_s$ and, respectively, $k_u$, so the corank of their tangency satisfies: $\mathcal C\le\min\{k_s,k_u\}\le\frac{k}{2}$. We will call the regions in the space of $C^r$-maps $(r=2,...,\infty,\omega)$ where tangencies of corank $\mathcal C>1$ are $C^2$-robust the \textit{Asaoka-Barrientos-Raibekas domain} (or, abbreviated, the \textit{$ABR$-domain}).

\begin{figure}[h]
\center{\includegraphics[width=0.5\linewidth]{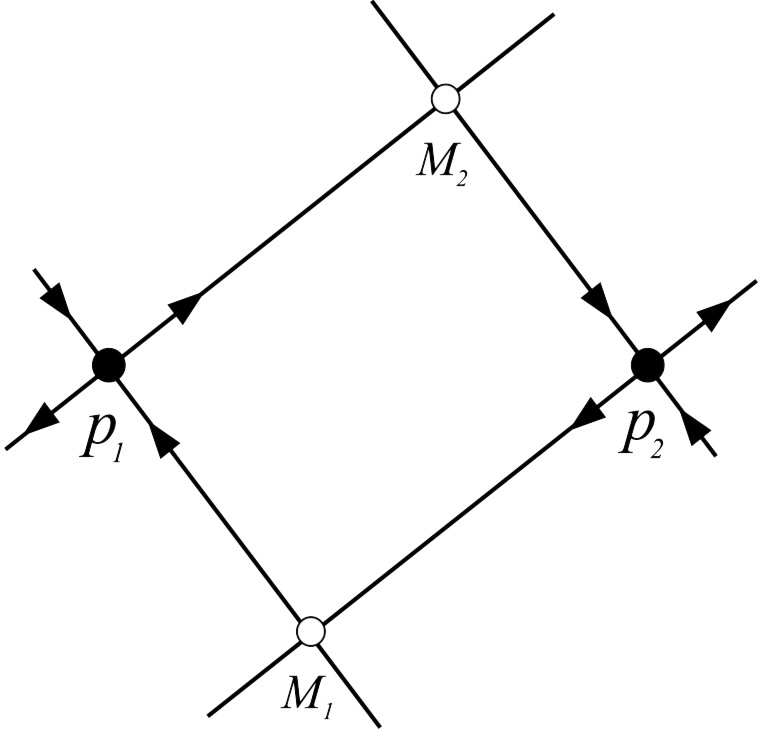}}
\caption{Homoclinically related non-trivial basic set $\Lambda$ and bi-focus periodic orbit $L$: $W^s(p_1)$ intersects $W^u(p_2)$ at the point $M_1$ and $W^u(p_1)$ intersects $W^s(p_2)$ at the point $M_2$, where $p_1\in\Lambda$ and $p_2\in L$.}
\label{Figure3}
\end{figure}

We consider a subdomain of the $ABR$-domain which is defined as follows. Let $f$ be a map from the $ABR$-domain which has a non-trivial basic set $\Lambda_f$ exhibiting a $C^2$-robust tangency of corank 2. Assume (see Fig. \ref{Figure3}) that the basic set $\Lambda_f$ is \textit{homoclinically related}\footnote{Recall that two basic sets $\Lambda_1$ and $\Lambda_2$ are said to be homoclinically related if there exist points $p_1\in\Lambda_1$ and $p_2\in\Lambda_2$ such that $W^s(p_1)\cap W^u(p_2)\not=\varnothing$ and $W^u(p_1)\cap W^s(p_2)\not=\varnothing$, and these intersections are transverse. } to a bi-focus periodic orbit $L_f$ (the homoclinic relation entails \cite{Smale1}, \cite{Shilnikov2} that there exist non-trivial basic sets which include $\Lambda_f$, $L_f$, and heteroclinic orbits connecting them). Maps with these properties obviously form open regions in the space of smooth systems. We will call such regions the \textit{$ABR^*$-domain}. 

\begin{theorem}\label{theorem2}
Given any $f\in ABR^*$ and any basic set $\Omega_f$ containing $\Lambda_f$ and $L_f$, arbitrarily $C^r$-close to $f$ there exists a map $g$ such that the continuation $\Omega_g$ has infinitely many orbits of corank-2 tangencies of every order between the stable and unstable manifolds of every pair of its periodic orbits.
\end{theorem}

The map $f$ belongs to the $ABR^*$-domain which is $C^r$-open, so Theorem \ref{theorem2} immediately implies Theorem \ref{theorem1}. Now let us explain how we prove Theorem \ref{theorem2} and why Theorem \ref{theorem3} is the key ingredient in the construction of the map $g$ (a detailed proof is given in Section \ref{subsection2.2}). Let $f\in ABR^*$ and let $\Lambda$ exhibit a $C^2$-robust corank-2 tangency $\Gamma_0$. We obtain the map $g$ as a result of applying a countable number of successive perturbations to the map $f$, each of which leaves a map in the $ABR^*$-domain and can be made arbitrarily $C^r$-small (so the total pertrubation is arbitrarily $C^r$-small).

\begin{enumerate}

\item The basic set $\Lambda$ and the periodic orbit $L$ are in the same basic set $\Omega$, therefore $W^s(L)$ and $W^u(L)$ accumulate to stable manifold and unstable manifold, respectively, of every point from $\Lambda$. It implies that adding a $C^r$-small perturbation which splits the initial tangency, we obtain an orbit of corank-2 homoclinic tangency between $W^s(L)$ and $W^u(L)$. The robustness (we are in the $ABR^*$-domain) entails that, in addition to the newly created corank-2 tangency between $W^s(L)$ and $W^u(L)$, we also have a corank-2 tangency between $W^s(\Lambda)$ and $W^u(\Lambda)$, splitting of which produces one more corank-2 tangency between $W^s(L)$ and $W^u(L)$. Repeating this procedure $h=2^{\frac{(N-1)(N+4)}{2}}$ times, where $N$ can by any natural number, we get $h$ orbits of corank-2 homoclinic tangency between $W^s(L)$ and $W^u(L)$ (along with some orbit of a robust corank-2 tangency between $W^s(\Lambda)$ and $W^u(\Lambda)$).

\item Applying Theorem \ref{theorem3}, we get an orbit of corank-2 homoclinic tangency of order $N$ between $W^s(L)$ and $W^u(L)$. 

\item Similarly to step 1, we add $C^r$-small perturbation to split the tangency obtained at the previous step. As a result, we get an orbit of corank-2 tangency of order $N$ between the stable and unstable manifolds of any given pair of periodic orbits of the set $\Omega$. 

\item The orbit obtained in step 3 coexists with some orbit of a robust corank-2 tangency between $W^s(\Lambda)$ and $W^u(\Lambda)$, so we can repeat steps 1-3 countably many times to obtain infinitely many orbits of corank-2 tangencies of every order between the stable and unstable manifolds of every pair of periodic orbits of the basic set $\Omega$, as claimed in Theorem \ref{theorem2}.

\end{enumerate}

If we exclude the real-analytic case from consideration, then the presence of high order tangencies of corank 2 in the system makes it possible (by adding an arbitrarily $C^r$-small perturbation to a map with a high order tangency of corank-2) to make the stable and unstable manifolds locally coincide along a two-dimensional disk. 

\begin{corollary}\label{corollary1}
Let $f$ be a $C^r$ map $(r=2,...,\infty)$ from the $ABR^*$-domain. Then arbitrarily $C^r$-close to $f$ there exists a map $h$ such that in a basic set $\Omega_h$ the intersection of the stable and unstable manifolds of every pair of periodic orbits contains a two-dimensional disk.
\end{corollary}

A similar result for corank-1 tangencies (the coincidence of stable and unstable manifolds along a curve) is obtained in \cite{arbitraryorder}, \cite{arbitraryorderconservative}, \cite{Kaloshin1}, \cite{Mints1}. It was used in \cite{Kaloshin1} for showing the genericity of the superexponential growth of the number of periodic points and in \cite{highsmoothness} for the construction of maps which cannot be topologically conjugate to any diffeomorphism of a higher smoothness.  

\subsection{Universal 2-dimensional dynamics}\label{subsection1.2}

Results from the previous section can be applied to the construction of maps with universal two-dimensional dynamics. By a 2-universal map one means that its iterations approximate the dynamics of every map of a two-dimensional disk arbitrarily well. In order to rigorously define the concept of universal dynamics, we follow the scheme from \cite{Turaev2} which provides a description for arbitrarily long iterations of a map on arbitrarily small spatial scales.

Let $\mathcal M^k$ be a smooth or analytic $k$-dimensional $(k\ge 2)$ manifold, let $F:\mathcal M^k\rightarrow \mathcal M^k$ be a $C^r$ diffeomorphism $(r=2,.\ldots,\infty,\omega)$. Let $U^k$ be the closed unit ball in $\mathbb R^k$, $B^k\subset \mathcal M^k$ be any small closed ball and $h: U^k\rightarrow B^k$ be a $C^r$ diffeomorphism. Given a positive $n$, the map $F^n|_{B^k}$ is a return map if $F^n(B^k)\cap B^k\not=\varnothing$. One can assume that the diffeomorphism $h: U^k\rightarrow B^k$ admits a $C^r$-smooth extension $H$ onto some larger ball $V^k\supset U^k$ such that $F^n(B^k)\subset H(V^k)$. By construction, the return map $F^n|_{B^k}$ is $C^r$-smoothly conjugate with the map $F_{n,H}=H^{-1}\circ F^n\circ H$. Thus, one obtains a $C^r$ map $F_{n,H}: U^k\rightarrow \mathbb R^k$ which is defined by the choice of the number of iterations $n$ and by the choice of the embedding $H$ (which also fixes the ball $B^k=H(U^k)$ on the manifold $\mathcal M^k$). The maps $F_{n,H}$ are called \textit{renormalized iterations} of $F$. The set $\bigcup\limits_{n,H} F_{n,H}$ of all possible renormalized iterations is called \textit{the dynamical conjugacy class} of $F$. 

\begin{definition}\label{definition1}

Let $F:\mathcal M^k\rightarrow\mathcal M^k$ be a $C^r$ diffeomorphism and let $m$ be a natural number such that $1\le m< k$. We say that $F$ has $C^r$-universal $m$-dimensional dynamics if the $C^r$-closure of its dynamical conjugacy class contains all maps of the unit ball $U^k$ into $\mathbb R^k$ of the following form:
\begin{equation*}
\begin{aligned} 
\left(\overline X_1,\ldots,\overline X_{k-m}\right)&= 0, \\
\left(\overline X_{k-m+1},\ldots,\overline X_k\right)&= \Phi \left(X_{k-m+1},\ldots,X_k\right),\\
\end{aligned}
\end{equation*}
where $\Phi$ is an arbitrary $C^r$ map of the unit ball $U^m$ into $\mathbb R^m$. 

\end{definition}

According to \cite{arbitraryorder},\cite{arbitraryorderconservative},\cite{Mints1}, for any $k$-dimensional $C^r$ diffeomorphism $(k\ge 2, r=2,...,\infty,\omega)$ with homoclinic tangency of corank 1, there exists $C^r$-small perturbation which produces a map with \textit{$C^r$-universal one-dimensional dynamics}. It entails that such universal maps are residual in the Newhouse regions of $k$-dimensional diffeomorphisms. For the $ABR^*$-domain we prove 
 
\begin{theorem}\label{theorem4}
In the $ABR^*$-domain $(r=2,\ldots,\infty,\omega)$ maps having $C^r$-universal two-dimensional dynamics form a residual subset\footnote{Let us recall that a set $A$ is said to be residual in a certain set $B$ if it is a countable intersection of open and dense subsets in $B$. In the real-analytic case $r=\omega$, residual set can be defined as follows: a set $A$ is residual in a certain subset $B$ of space of real-analytic maps, if given any map $f$ from $B$ and any compact subset $K\subset\mathcal M^k$, there exists a complex neighborhood $Q$ of $K$ such that the intersection of $A$ with some open neighborhood $X$ of $f$ in space of maps holomorphic on $Q$ is the intersection of a countable collection of open and dense subsets of $X$.}. 
\end{theorem}

Note that maps $\Phi$ in Definition \ref{definition1} do not need to be invertible, i.e., the variety of dynamics displayed by the 2-universal maps in the sense of this Definition is greater than for the 2-universal maps which are generic in the so-called absolute Newhouse domain (see \cite{Turaevintro1}, \cite{Turaevintro2}) -- for the universal maps constructed in \cite{Turaevintro2}, the dynamical conjugacy class is only dense among all orientation preserving diffeomorphisms.  

The concept of universal dynamics and Theorem \ref{theorem4} can be generalized to finite-parameter families of smooth and real-analytic maps. Let $\mathcal P_{l,r}$ be the space of all $l$-parameter families $(l\in\mathbb N)$ of diffeomorphisms of the manifold $\mathcal M^k$, of class $C^r$ $(r=2,.\ldots,\infty,\omega)$ with respect to the coordinates and the parameters $\varepsilon=(\varepsilon_1,\ldots,\varepsilon_l)$. Let $\mathcal D^l\subset \mathbb R^l$ be the space where the parameters $\varepsilon$ are defined. Then any family of maps $F_{\varepsilon}\in \mathcal P_{l,r}$ can be considered as a map from $\mathcal M^k\times \mathcal D^l$ into $\mathcal M^k$ which acts as $F_{\varepsilon}$ on $\mathcal M^k$ and as the identity map on $\mathcal D^l$. Making use of this representation, in exactly the same way as it was done for maps, one can define \textit{renormalized iterations} and \textit{the dynamical conjugacy class} for any family $F_{\varepsilon}\in\mathcal P_{l,r}$.
 
\begin{definition}\label{definition1parametricfamily}

Let a parametric family of maps $F_{\varepsilon}$ belongs to $\mathcal P_{l,r}$  and let $m$ be a natural number such that $1\le m< k$. We say that $F_{\varepsilon}$ has $C^r$-universal $m$-dimensional dynamics if the $C^r$-closure of its dynamical conjugacy class contains all $l$-parameter families of maps of the unit ball $U^k$ into $\mathbb R^k$ of the following form:
\begin{equation*}
\begin{aligned} 
\left(\overline X_1,\ldots,\overline X_{k-m}\right)&= 0, \\
\left(\overline X_{k-m+1},\ldots,\overline X_k\right)&= \Phi \left(X_{k-m+1},\ldots,X_k,\varepsilon_1,\ldots,\varepsilon_l\right),\\
\end{aligned}
\end{equation*}
where $\Phi$ is an arbitrary $C^r$ map of the unit ball $U^{m+l}$ into $\mathbb R^m$. 

\end{definition}

Let us denote by $\mathcal P_{l,r}(ABR^*)$ the subspace of the space $\mathcal P_{l,r}$ such that all maps of each family $F_{\varepsilon}\in \mathcal P_{l,r}(ABR^*)$ belong to the $ABR^*$-domain. Then, absolutely similarly to the proof of Theorem \ref{theorem4}, one can prove

\begin{theorem}\label{theorem5}
In the space $\mathcal P_{l,r}(ABR^*)$ $(l\in\mathbb N, r=2,\ldots,\infty,\omega)$ parameteric families having $C^r$-universal two-dimensional dynamics form a residual subset.
\end{theorem}

\subsection{Further discussion}\label{subsection1.3}

One can characterize the Newhouse domain as an open domain in the space of $C^r$-systems $(r=2,\ldots,\infty,\omega)$ where each system has a non-trivial basic set which exhibits a $C^2$-robust tangency of corank 1 (also called \textit{wild hyperbolic set}, see \cite{Newhouse3}). If the dimension of the ambient manifold $k\ge 4$, then we can select an open subdomain of the Newhouse domain where each system has wild hyperbolic set containing bi-focus periodic orbit. We will call such subdomain a \textit{bi-focus Newhouse domain}. 

In the paper in preparation \cite{Mints1}, we show that a $C^r$-small perturbation of a map with a bi-focus periodic orbit, whose stable and unstable manifolds contain a corank-1 homoclinic tangency, leads to the creation of an infinite number of orbits of corank-2 homoclinic tangency, i.e. corank-2 tangencies are persistent and turn out to be very natural phenomenon. Combining these results with an algorithm provided by Theorem \ref{theorem3} gives the following statement.\\ 

\textit{In the space of $k$-dimensional $C^r$ maps, where $k\ge 4$ and $r=2,...,\infty,\omega$, in any neighborhood of a map such that it has a bi-focus periodic orbit whose stable and unstable manifolds are tangent, there exist bi-focus Newhouse regions in which}
 \begin{enumerate}
 
\item \textit{maps with infinitely many orbits of corank-2 homoclinic tangencies of every order form a dense subset;}
 
\item \textit{maps having universal two-dimensional dynamics form a residual subset.}
 
\end{enumerate}
 
\section{High order tangencies of corank 2}\label{section2}

In this section, we prove Theorem \ref{theorem2} which immediately implies our main result, Theorem \ref{theorem1}. We assume that Theorem \ref{theorem3} holds, the proof of which is independent and given in Section \ref{subsection3.4}. 

In Section \ref{subsection2.1}, we make some preliminary constructions and prove Lemma \ref{lemmasplittingofcorank2tangency} which makes it possible to \textit{split generically and independently} any finite number of corank-2 tangencies. In Section \ref{subsection2.2}, we prove Lemma \ref{lemmamanytangencies} which shows that generic system from $ABR^*$-domain have an \textit{abundance of corank-2 homoclinic tangencies to the bi-focus periodic orbit}. Based on Lemma \ref{lemmamanytangencies} and Theorem \ref{theorem3}, we prove Theorem \ref{theorem2} and Corollary \ref{corollary1}.

\subsection{Splitting of corank-2 tangencies}\label{subsection2.1}

Let $f$ be a $C^r$ map $(r=2,\ldots,\infty,\omega)$ with a basic set $\Lambda$ of type $(k_s,k_u)$. Let $p,q$ be two points (possibly coinciding) such that $p,q\in\Lambda$ and $W^s(p), W^u(q)$ form a corank-2 tangency $\Gamma$ at the point $M$. Then from the definition of a corank-2 tangency (see Section \ref{subsection1.1}) it follows that near $M$ there are coordinates $(x_1,\ldots,x_{k_s},y_1,\ldots,y_{k_u})$ in which $W^s(p)$ and $W^u(q)$ locally have the form
\begin{equation}\label{corank2tangency}
\begin{aligned} 
&W^s(p)=\{y_1=\ldots=y_{k_u}=0\},\\
&W^u(q)=\{(y_1,y_2)=G(x_1,x_2),x_3=\ldots=x_{k_s}=0\},
\end{aligned}
\end{equation}
where $G$ is a $C^r$ map such that $G(0)=0$ and $DG(0)=0$ (it is the generalization of the form for corank-1 tangencies, see the discussion in \cite{bifurcationsandstability}, section II.6). When $W^s(p)$ and $W^u(q)$ are given in the form \eqref{corank2tangency}, we can define an order of corank-2 tangency as follows.

\begin{definition}\label{definition2}
We say that the corank-2 tangency $\Gamma$ is of order $n<r$, $n\in\mathbb N$, if 
\begin{enumerate}

\item $\frac{\partial^{i+j} G(0)}{\partial x^i_1\partial x^j_2}=0$ for all $i,j\in\{0\}\cup\mathbb N$ such that $(i+j)\in\{1,...,n\}$;

\item $\frac{\partial^{n+1} G(0)}{\partial x^i_1\partial x^j_2}\not=0$ for at least one pair of $i,j\in\{0\}\cup\mathbb N$ such that $i+j=n+1$.

\end{enumerate}

If $r$ is finite and the first condition holds for $n=r$, or $r=\infty$ and the first condition holds for all $n\in\mathbb N$, then we say that the tangency $\Gamma$ is $C^r$-flat. 

\end{definition}

Assume that the corank-2 tangency $\Gamma$ is of order $n$ and consider $C^r$-small perturbations of the map $f$. For all such perturbations, in suitable coordinates, the equations for $W^s(p)$ and $W^u(q)$ again can be written in the form \eqref{corank2tangency} with the map $G$ given by
\begin{equation}\label{mapFperturbations}
\begin{aligned} 
&y_1 = \sum\limits_{j=0}^n\sum\limits_{i=0}^{j} \mu_{j,i}x_1^{j-i} x_2^{i}+\sum\limits_{i=0}^{n+1}A_i x_1^{n+1-i}x_2^{i}+o\left((x_1^2+x_2^2)^{\frac{n+1}{2}}\right), \\
&y_2=\sum\limits_{j=0}^n\sum\limits_{i=0}^{j} \nu_{j,i} x_1^{j-i} x_2^{i}+\sum\limits_{i=0}^{n+1}B_i x_1^{n+1-i}x_2^{i}+o\left((x_1^2+x_2^2)^{\frac{n+1}{2}}\right), \\
\end{aligned}
\end{equation}
where all the coefficients depend continuously in the $C^r$-topology, $\mu_{j,i}$ and $\nu_{j,i}$ are smooth functionals which are small in absolute value, and $\sum\limits_{i=0}^{n+1} |A_i|+\sum\limits_{i=0}^{n+1} |B_i|>0$. 

Assume that the map $f$ is included in an arbitrary finite-parameter family of $C^r$ maps $f_{\varepsilon}$ $(f_0=f)$. Let us fix a system of coordinates such that for all small $\varepsilon$ the map $G$ is given by \eqref{mapFperturbations}. Denote by $\mu$, $\nu$ the vectors of coefficients $\mu_{j,i}$, $\nu_{j,i}$, respectively, by $\chi(\mu,\nu)=n^2+3n+2$ the total number of coefficients $\mu_{j,i}$, $\nu_{j,i}$. 

\begin{definition}\label{definition3}
Corank-2 tangency of order $n$ is said to be split generically if at all small $\varepsilon$ the following condition holds:
\begin{equation}\label{corank2splitgenerically} 
\text{rank}\;\frac{\partial(\mu,\nu)}{\partial\varepsilon}=\chi(\mu,\nu). 
\end{equation}
\end{definition}

Now assume that the map $f$ has corank-2 tangencies $\Gamma_1,\ldots,\Gamma_h$ of orders $n_1,\ldots n_h$, respectively, and $f$ is included in an arbitrary finite-parameter family of $C^r$ maps $f_{\varepsilon}$ $(f_0=f)$. As above, for each tangency $\Gamma_l$ $(l=1,\ldots,h)$ we define smooth functionals $\mu^l_{j,i}$, $\nu^l_{j,i}$ $(l=1,\ldots,h)$ which determine its splitting. We denote by $\mu^l$, $\nu^l$ the vectors of $\mu^l_{j,i}$, $\nu^l_{j,i}$, respectively, by $\chi^l(\mu^l,\nu^l)$ the total number of $\mu^l_{j,i}$, $\nu^l_{j,i}$.  

\begin{definition}\label{definition3.5}
Corank-2 tangencies $\Gamma_1,\ldots,\Gamma_h$ are said to be split independently if at all small $\varepsilon$ the following condition holds:
\begin{equation}\label{corank2splitindependently}
\text{rank}\;\frac{\partial(\mu^1,\nu^1,\ldots,\mu^h,\nu^h)}{\partial(\varepsilon_1,\ldots,\varepsilon_h)}=\sum\limits_{l=1}^h\chi_l(\mu^l,\nu^l).
\end{equation} 
\end{definition}

The \textit{independence of the splitting} means that, for each of the tangencies $\Gamma_1,\ldots,\Gamma_h$ under consideration, there exists a smooth manifold $S_i(\Gamma_i)$ in space of parameters $\varepsilon$ such that when $\varepsilon$ varies within $S_i(\Gamma_i)$ the tangency $\Gamma_i$ splits generically, while the other tangencies are not split at all, and the manifolds $S_i(\Gamma_i)$ and $S_j(\Gamma_j)$ corresponding to different tangencies intersect at $\varepsilon=0$ transversely.

Further, we construct perturbations that allow to \textit{split generically} and \textit{independently} any finite number of corank-2 tangencies.

\begin{lemma}\label{lemmasplittingofcorank2tangency}
Let a $C^r$ map $f$, where $r=2,...,\infty,\omega$, have corank-2 tangencies $\Gamma_1,\ldots,\Gamma_h$. Then there exists an analytic finite-parameter family of real analytic maps $\mathcal F_{\varepsilon}$ such that in the family $\mathcal F_{\varepsilon} \circ f$ the tangencies $\Gamma_1,...,\Gamma_h$ split generically and independently.
\end{lemma}

\begin{demo}

\underline{Construction of $C^{\infty}$ map $\mathcal F_{\varepsilon}$.}

First we prove the lemma for $h=1$. Let the corank-2 tangency $\Gamma_1$ be of order $n_1$. Then stable and unstable manifolds forming $\Gamma_1$ are given by \eqref{corank2tangency}, the map $G$ has the form 
\begin{equation}\label{lemmasplittingofcorank2tangencyeq1}
\begin{aligned} 
&y_1 = \sum\limits_{i=0}^{n+1}A_i x_1^{n+1-i}x_2^{i}+o\left((x_1^2+x_2^2)^{\frac{n+1}{2}}\right), \\
&y_2=\sum\limits_{i=0}^{n+1}B_i x_1^{n+1-i}x_2^{i}+o\left((x_1^2+x_2^2)^{\frac{n+1}{2}}\right), \\
\end{aligned}
\end{equation}
where $\sum\limits_{i=0}^{n+1} |A_i|+\sum\limits_{i=0}^{n+1} |B_i|>0$.

Fix a small $\rho>0$ and define a $C^{\infty}$-smooth cut-off function $\theta_{\rho}=\theta_{\rho}(x_1,\ldots,x_{k_s},y_1,\ldots,y_{k_u})$ such that it is equal to 1 inside the $\rho$-ball $B^1_{M}$ and vanishes identically outside the $2\rho$-ball $B^2_{M}$ (these balls are centered at the point $M$). Let $F_{\varepsilon_1}$ be a $C^{\infty}$ map depending on $\chi_1(\mu^1,\nu^1)$ parameters and acting as  
\begin{equation}\label{lemmasplittingofcorank2tangencyeq2}
\begin{aligned} 
&(x_1,\ldots,x_{k_s},y_1,\ldots,y_{k_u})\rightarrow (x_1,\ldots,x_{k_s},y_1+\theta{\rho}\cdot \sum\limits_{j=0}^n\sum\limits_{i=0}^{j} \mu^1_{j,i} x_1^{j-i} x_2^{i},\\
&y_2+\theta{\rho}\cdot \sum\limits_{j=0}^n\sum\limits_{i=0}^{j} \nu^1_{j,i} x_1^{j-i} x_2^{i},\ldots,y_{k_u}).
\end{aligned} 
\end{equation} 
By construction, the map $F_{\varepsilon_1}$ is equal to the identity outside a small neighborhood of the point $M$ for all $\varepsilon_1$; at $\varepsilon_1=0$ it is equal to the identity everywhere. 

Consider the family $F_{\varepsilon_1}\circ f$. In this family the equation for $W^s(p_1)$ does not change while the equation for $W^u(p_2)$ takes the form $W^u(p_2)=\{(y_1,y_2)=G^{new}(x_1,x_2),x_3=\ldots=x_{k_s}=0\}$, where $G^{new}(x_1,x_2)$ coincides with \eqref{mapFperturbations} in a $\rho$-neighborhood of the point $M$. It is easily seen that for the family $F_{\varepsilon_1}\circ f$ the condition \eqref{corank2splitgenerically} holds (the corresponding $\chi_1\times \chi_1$-matrix is the identity matrix), therefore the tangency $\Gamma_1$ is split generically in the family $F_{\varepsilon_1}\circ f$.

Now we prove the lemma for an arbitrary finite $h$. Let the corank-2 tangencies $\Gamma_1,\ldots,\Gamma_h$ have orders $n_1,\ldots,n_h$. For each tangency $\Gamma_l$ $(l=1,\ldots,h)$ we define smooth functionals $\mu^l_{j,i}$, $\nu^l_{j,i}$ $(l=1,\ldots,h)$ which determine its splitting. Consider a $C^{\infty}$ map $\mathcal F_{\varepsilon}$ given by the formula:
\begin{equation*}
\mathcal F_{\varepsilon}=F_{\varepsilon_1}+\ldots+F_{\varepsilon_h},
\end{equation*} 
where the maps $F_{\varepsilon_1},\ldots,F_{\varepsilon_h}$ correspond to $\Gamma_1,\ldots,\Gamma_h$ and are localized in sufficiently small neighborhoods of the given tangencies. Now consider a finite-parameter family $\mathcal F_{\varepsilon}\circ f$. The previous reasoning implies that for such family the equality \eqref{corank2splitindependently} holds. Thus, in the family $\mathcal F_{\varepsilon}\circ f$ all the tangencies $\Gamma_1,\ldots,\Gamma_h$ are split generically and independently. 

\underline{Construction of real-analytic map $\mathcal F_{\varepsilon}$.}

Let a map $\mathcal F^{new}_{\varepsilon}$ be suffiently close (in the $C^r$ metric) analytic approximation of the map $\mathcal F_{\varepsilon}$. Then the equality \eqref{corank2splitindependently} is preserved for the family $\mathcal F^{new}_{\varepsilon}\circ f$. It completes the proof of the lemma.

\end{demo}

\subsection{Abundance of high order tangencies of corank 2 in the $ABR^*$-domain}\label{subsection2.2}

\begin{lemma}\label{lemmamanytangencies}
Let a real-analytic map $f_0\in ABR^*$, that is, $f_0$ has homoclinically related bi-focus periodic orbit $L_{f_0}$ and non-trivial basic set $\Lambda_{f_0}$ exhibiting $C^2$-robust tangency of corank 2. Then there exists an analytic finite-parameter family of real-analytic maps $f_{\varepsilon}$ and a sequence $\varepsilon_i\xrightarrow[i\rightarrow+\infty]{}0$ such that $W^s(L_{f_{\varepsilon_i}})$ and $W^u(L_{f_{\varepsilon_i}})$ form $h$ homoclinic tangencies $\Gamma_1,\ldots,\Gamma_h$ of corank 2, where $h$ can be chosen to be any natural number. 

In addition, perturbation can be chosen in such a way that it does not affect any given finite number of other tangencies in the system. 
\end{lemma}

\begin{demo}

Since $\Lambda_{f_0}$ and $L_{f_0}$ are homoclinically related, then there exists a basic set $\Omega_{f_0}$ containing them. It entails that $W^s(L_{f_0})$ and $W^u(L_{f_0})$ accumulate to stable manifold and unstable manifold, respectively, of every point from $\Lambda_{f_0}$. Moreover, this property holds for all $C^{\omega}$-small perturbations of the map $f_0$.

Let $\Lambda_{f_0}$ have a corank-2 tangency $\Gamma_0$ between $W^s(p_0)$ and $W^u(q_0)$, where points $p_0,q_0\in\Lambda_{f_0}$. By Lemma \ref{lemmasplittingofcorank2tangency}, we can include the map $f_0$ in an analytic finite-parameter family of real-analytic maps in which the tangency $\Gamma_0$ splits generically and any given finite number of other tangencies is not affected. Since invariant manifolds $W^s(L_{f_0})$ and $W^u(L_{f_0})$ accumulate to invariant manifolds $W^s(p_0)$ and $W^u(q_0)$, respectively, then in this family, arbitrarily close to $f_0$, we find a map $f_1$ which has a corank-2 homoclinic tangency $\Gamma_1$ between $W^s(L_{f_1})$ and $W^u(L_{f_1})$. Also, since the map $f_1\in ABR^*$-domain, it has a corank-2 tangency $\Gamma_{00}$ between $W^s(p_1)$ and $W^u(q_1)$, where points $p_1,q_1\in\Lambda_{f_1}$. 

Again, making use of Lemma \ref{lemmasplittingofcorank2tangency}, we include the map $f_1$ in an analytic finite-parameter family of real-analytic maps in which the tangency $\Gamma_{00}$ splits generically, and the tangency $\Gamma_1$ along with any given finite number of other tangencies is not affected. In this family, arbitrarily close to $f_1$, we find a map $f_2$ which has corank-2 homoclinic tangencies $\Gamma_1,\Gamma_2$ between $W^s(L_{f_2})$ and $W^u(L_{f_2})$ and a corank-2 tangency $\Gamma_{000}$ between $W^s(p_2)$ and $W^u(q_2)$, where points $p_2,q_2\in\Lambda_{f_2}$.
 
Performing this procedure $h$ times, we obtain a map $f_h$ which has $h$ homoclinic tangencies $\Gamma_1,\ldots,\Gamma_h$ of corank 2 between $W^s(L_{f_h})$ and $W^u(L_{f_h})$. By construction, the whole perturbation preserves any given finite number of other tangencies for the map $f_0$ and can be chosen arbitrarily small in the $C^{\omega}$ topology.

\end{demo}

\begin{demo3}

Let a $C^r$ map $f\in ABR^*$-domain $(r=2,\ldots,\infty,\omega)$, that is, $f$ has homoclinically related bi-focus periodic orbit $L$ and non-trivial basic set $\Lambda$ exhibiting $C^2$-robust tangency of corank 2. Let $U$ be a neighborhood of $\mathcal L\subset W^s(\Lambda)\cup W^u(\Lambda)$ which contains orbits of robust tangencies of stable and unstable manifolds of the set $\Lambda$ for all maps sufficiently close to $f$ (see definition of robust tangency in Section \ref{subsubsection1.1.2}). Since the neighborhood $U$ can be chosen arbitrarily small, we choose it in such a way that it does not contain the orbit $L$. We will consider that the map $f$ is real-analytic (if $f$ is only smooth, then arbitrarily close to $f$ one can always find a map $f^{new}$ which is real-analytic). 

Let $\Omega$ be some basic set containing $L$ and $\Lambda$. Let $n_1,n_2,\ldots$ be an arbitrary infinite sequence of natural numbers, and let $(M_{1,1};M_{1,2}),(M_{2,1};M_{2,2}),(M_{3,1};M_{3,2}),\ldots$ be an arbitrary sequence of pairs of periodic points from the basic set $\Omega$. We will show that there exists a perturbation of $f$ leading to the map with an infinite sequence of homoclinic/heteroclinic tangencies of corank 2, and these will be exactly the tangencies of orders $n_k$ between manifolds $W^s(M_{k,1})$ and $W^u(M_{k,2})$, where $k\in\mathbb N$. This perturbation can be as small as we need (in the $C^{\omega}$ topology). The existence of such perturbation immediately provides a proof of the theorem (one should allow the numbers $n_k$ take all natural values infinitely many times).

Take an arbitrarily small $\delta>0$, and let $\delta_k>0$ be such that $\delta_1+\delta_2+\ldots=\delta$. We will construct a sequence of maps $f_k$, where $f_0\equiv f$, such that $f_k$ has one orbit of $C^2$-robust tangency of corank 2 between $W^s(p_k)$, $W^u(q_k)$ for some points $p_k,q_k\in\Lambda$ and $k$ orbits of homoclinic/heteroclinic tangency of corank 2: one orbit of tangency of order $n_1$ between $W^s(M_{1,1})$ and $W^u(M_{1,2})$, one orbit of tangency of order $n_2$ between $W^s(M_{2,1})$ and $W^u(M_{2,2})$, etc. The maps $f_k$ will have the property that the distance between $f_{k+1}$ and $f_k$ will be less than $\delta_k$. The sequence of maps $f_k$ will have a limit $f^*$ which is at a distance less than $\delta$ from $f$, and possess the required sequence of tangencies.

Thus, to prove the theorem, we have to show that given a map $f_k$ with $C^2$-robust tangency of corank 2 and $k\ge 0$ orbits of homoclinic/heteroclinic tangency of corank 2, one can perturb it and get the map $f_{k+1}$ with the following properties:

\begin{enumerate}

\item $f_{k+1}$ has $C^2$-robust tangency of corank 2 between $W^s_{loc,n}(p_{k+1})$, $W^u_{loc,n}(q_{k+1})$ for some points $p_{k+1},q_{k+1}\in\Lambda$;

\item $f_{k+1}$ preserves $k$ tangencies of corank 2 along with their orders;

\item $f_{k+1}$ has one more corank-2 tangency $\Gamma$ between $W^s(M_{k+1,1})$ and $W^u(M_{k+1,2})$ of the given order $N=n_{k+1}$.

\end{enumerate}

We will construct such a perturbation of the map $f_k$ as a finite sequence of perturbations each of which can be made arbitrarily $C^{\omega}$-small (since all these perturbations are given by Lemma \ref{lemmamanytangencies} and Theorem \ref{theorem3} which is based on Lemmas \ref{lemmahsotangency}, \ref{lemmahotangency}, \ref{lemmahotangency1}). Therefore the resulting perturbation will be of class $C^{\omega}$ and of total size less than $\delta_{k+1}$. 

By Lemma \ref{lemmamanytangencies}, perturbing the map $f_k$ we can create another map $\overline f_k$ such that it has $C^2$-robust tangency of corank 2 between $W^s(p_{k+1})$, $W^u(q_{k+1})$ for some points $p_{k+1},q_{k+1}\in\Lambda$ (the perturbation does not lead out of the $ABR^*$-domain), and $h=2^{\frac{(N-1)(N+4)}{2}}$ homoclinic tangencies $\Gamma_1,\ldots,\Gamma_h$ of corank 2 between $W^s(L)$ and $W^u(L)$. In addition, this perturbation does not affect any given finite number of tangencies in the system. Applying Theorem \ref{theorem3} to the map $\overline f_k$, we get a map $\tilde f_k$ with a corank-2 homoclinic tangency $\tilde\Gamma$ of order $N$ between $W^s(L)$ and $W^u(L)$, while also not affecting any given finite number of tangencies in the system.

Let us recall that the points $L, M_{k+1,1}$ and $M_{k+1,2}$ belong to the same non-trivial basic set $\Omega$, therefore invariant manifolds $W^s(M_{k+1,1})$ and $W^u(M_{k+1,2})$ accumulate to invariant manifolds  $W^s(L)$ and $W^u(L)$, respectively. Adding arbitrarily $C^{\omega}$-small local perturbation which splits the tangency $\tilde\Gamma$ generically (such perturbation is given by Lemma \ref{lemmasplittingofcorank2tangency}), we get a map $f_{k+1}$ with the sought corank-2 tangency $\Gamma$ of order $N$ between $W^s(M_{k+1,1})$ and $W^u(M_{k+1,2})$. 

Let us emphasize that if $M_{k+1,1}$ and $M_{k+1,2}$ belong to the set $\Lambda$, then at least one point of the orbit of tangency $\Gamma$ does not belong to the set $U$. Indeed, by construction, at least one point of the orbit of tangency $\Gamma$ would be close to $L$ (proximity to $L$ can be controlled by the size of the perturbation which splits the tangency $\tilde\Gamma$) which does not belong to $U$. Thus, the tangency $\Gamma$ cannot be $C^2$-robust therefore the above procedure of obtaining maps $f_k$ can be performed further while maintaining the tangency $\Gamma$.

This completes the proof of Theorem \ref{theorem2}. 

\end{demo3}

\begin{democor1}

Let $r=2,\ldots,\infty$. It directly follows from the proof of Theorem \ref{theorem2} that arbitrarily $C^r$-close to any $f\in ABR^*$ one can construct a $C^r$ map $\hat g$ with a basic set $\Omega_{\hat g}$ such that the stable and unstable manifolds of every pair of periodic orbits of $\Omega_{\hat g}$ contain an orbit of corank-2 tangency which is

\begin{enumerate}

\item $C^r$-flat if $r$ is finite;

\item of order $n$ if $r=\infty$, where $n$ can be as large as one needs.

\end{enumerate}

Each such tangency can be transformed to a two-dimensional disk by a perturbation which is arbitrarily small (in the $C^r$ topology) and which is localized in an arbitrarily small given neighborhood of one point of the orbit of tangency under consideration. In addition, one can localize these perturbations in such a way that they do not affect any other considered tangencies. As a result of sequential application of such perturbations, arbitrarily $C^r$-close to $\hat g$, one obtains a map $h\in ABR^*$ such that it has non-trivial basic set $\Omega_{h}$ in which the intersection of the stable and unstable manifolds of every pair of periodic orbits contains a two-dimensional disk.

Thus, near any map $f\in ABR^*$ one obtains the map $h\in ABR^*$ with the required property. It completes the proof.
\end{democor1}

\section{High order tangencies of corank-2 to bi-focus periodic orbit}\label{section3}

In this section, we prove Theorem \ref{theorem3}. Let us emphasize that all the results in this section are independent of the $ABR^*$-domain.

In our proofs, we "split" the dynamics of maps with homoclinic orbits into two parts, considering the so-called \textit{local} and \textit{global maps}. In Sections \ref{subsection3.1} and \ref{subsection3.2}, we provide explicit formulas for these maps. In Section \ref{subsection3.3}, we introduce \textit{index} of corank-2 tangency, which plays a leading role in the proof of Theorem \ref{theorem3}, and generalize ideas from Section \ref{subsection2.1}. In Section \ref{subsection3.3}, we provide key Lemmas \ref{lemmahsotangency}, \ref{lemmahotangency}, \ref{lemmahotangency1} which serve as a basis for obtaining high-order tangencies of corank 2. After it, we prove Theorem \ref{theorem3}. 

\subsection{Local map}\label{subsection3.1}

Let $f$ be a $C^r$ map $(r=2,\ldots,\infty,\omega)$ with a bi-focus periodic orbit $L=\{O,f(O),...,f^{p-1}(O)\}$ of period $p$. Denote by $\lambda_1,...,\lambda_{k_1}$, $\gamma_1,\ldots,\gamma_{k_2}$ the multipliers of $L$, and order them so that $|\gamma_{k_2}|\ge...\ge|\gamma_1|>1>|\lambda_1|\ge\ldots\ge|\lambda_{k_1}|$. Denote $\lambda_{1,2}=\lambda e^{\pm i\varphi}, \; \gamma_{1,2}=\gamma e^{\pm i\psi}$, where $\varphi,\psi\not=0,\pi$, and $\lambda=|\lambda_{1,2}|$, $\gamma=|\gamma_{1,2}|$. Note that the definition of a bi-focus periodic orbit implies that if $L$ has non-leading multipliers (i.e. the dimension of the phase space is more than 4), then
\begin{equation}\label{eqmult1}
\lambda>|\lambda_3|, \;\;\;\; \gamma<|\gamma_3|.
\end{equation}
We will assume that $\lambda\gamma\not=1$ since this can always be achieved by adding a $C^r$-small perturbation in a small neighborhood of $L$. We will consider the case $\lambda\gamma<1$ (if $\lambda\gamma>1$, then it is sufficient to consider a map $f^{-1}$ instead of $f$). 

Let $T_0$ be the restriction of the Poincar{\'e} map $f^p$ onto a small neighborhood $U$ of $O$. The map $T_0$ is said to be the \textit{local map}. We will consider such $C^r$ coordinates in the neighborhood $U$ that the fixed point $O$ of the local map $T_0$ is at the origin for all small $\varepsilon$. We will assume that the local stable and unstable manifolds of the point $O$ are straightened. Taking into account \eqref{eqmult1}, we can write the map $T_0$ in the form
\begin{equation}\label{eqmult2}
\begin{aligned} 
&\overline x_1=\lambda\left(x_1\cos(k\varphi)-x_2\sin(k\varphi)\right)+\ldots, \;\;\;\;\;\;\;\;\;\;\;\;\;\;\;\;\;\;\;\; \overline y_1=\gamma\left(y_1\cos(k\psi)-y_2\sin(k\psi)\right)+\ldots,\\
&\overline x_2=\lambda\left(x_1\sin(k\varphi)+x_2\cos(k\varphi)\right)+\ldots, \;\;\;\;\;\;\;\;\;\;\;\;\;\;\;\;\;\;\;\; \overline y_2=\gamma\left(y_1\sin(k\psi)+y_2\cos(k\psi)\right)+\ldots,\\
&\overline u=A u+\ldots,\;\;\;\;\;\;\;\;\;\;\;\;\;\;\;\;\;\;\;\;\;\;\;\;\;\;\;\;\;\;\;\;\;\;\;\;\;\;\;\;\;\;\;\;\;\;\;\;\;\;\;\;\;\;\;\;\;\;\; \overline v=B v+\ldots,
\end{aligned} 
\end{equation}
where the dots denote nonlinear terms; $u\in \mathbb R^{k_1-2}, v\in \mathbb R^{k_2-2}$; the eigenvalues of the matrices $A$ and $B$ are, respectively, the stable and unstable non-leading multipliers of the orbit $L$. 

We will consider the local map $T_0^k$ in the so-called \textit{main normal form}, which can always be reduced to by $C^r$ transformations of coordinates (see \cite{Gonchenko3}). Basically, this means the absence of some nonlinear terms in the Taylor expansion, which make the iterations of the local map too far from its linearization. Also, it allows us to rewrite the map $T_0^k$ in the \textit{cross form} in which all our proofs are carried out. Let us formulate the above as 

\begin{proposition}\label{proposition1} (\cite{Gonchenko3}, Lemma 7)
Let a $C^r$ map $f$, where $r=2,...,\infty,\omega$, have a bi-focus periodic orbit $L$. Then for all small $\varepsilon$ there exist $C^r$ coordinates in $U$ such that in these coordinates the map $T_0^k$ for all large $k$ has the following form
\begin{equation}\label{eqlocal3}
\begin{aligned} 
&x_{k1}=\lambda^k(\varepsilon)\cdot\left(x_{01}\cos(k\varphi(\varepsilon))-x_{02}\sin(k\varphi(\varepsilon))\right)+\hat\lambda^k\cdot\xi^1_k(x_{01},x_{02},y_{k1},y_{k2},u_0,v_k,\varepsilon),\\
&x_{k2}=\lambda^k(\varepsilon)\cdot\left(x_{01}\sin(k\varphi(\varepsilon))+x_{02}\cos(k\varphi(\varepsilon))\right)+\hat\lambda^k\cdot\xi^2_k(x_{01},x_{02},y_{k1},y_{k2},u_0,v_k,\varepsilon),\\
&y_{01}=\gamma^{-k}(\varepsilon)\cdot\left(y_{k1}\cos(k\psi(\varepsilon))+y_{k2}\sin(k\psi(\varepsilon))\right)+\hat\gamma^{-k}\cdot\eta^1_k(x_{01},x_{02},y_{k1},y_{k2},u_0,v_k,\varepsilon),\\
&y_{02}=\gamma^{-k}(\varepsilon)\cdot\left(-y_{k1}\sin(k\psi(\varepsilon))+y_{k2}\cos(k\psi(\varepsilon))\right)+\hat\gamma^{-k}\cdot\eta^2_k(x_{01},x_{02},y_{k1},y_{k2},u_0,v_k,\varepsilon),\\
&u_k=\hat\lambda^k\cdot\hat\xi_k(x_{01},x_{02},y_{k1},y_{k2},u_0,v_k,\varepsilon),\\
&v_0=\hat\gamma^{-k}\cdot\hat\eta_k(x_{01},x_{02},y_{k1},y_{k2},u_0,v_k,\varepsilon),
\end{aligned}
\end{equation}
where $\hat\lambda$ and $\hat\gamma$ are some constants such that $0<\hat\lambda<\lambda$, $\hat\gamma>\gamma$, and the functions $\xi^1_k,\xi^2_k, \eta^1_k,\eta^2_k, \hat\xi_k, \hat\eta_k$ are uniformly bounded for all $k$, along with the derivatives up to the order $r-2$.

\end{proposition}

\begin{remark}\label{remark1}
The proof from \cite{Gonchenko3} (Lemma 7) directly implies that the constants $\hat\lambda$ and $\hat\gamma$ can be chosen such that
\begin{equation}
  \label{eqlambdagamma}
   \hat\lambda<\lambda<\hat\gamma^{-1}<\gamma^{-1} \quad\mbox{and}\quad \gamma^{-2}<\hat\gamma^{-1}.
\end{equation}
\end{remark}

\subsection{Global map}\label{subsection3.2}

Now assume that the stable and unstable manifolds of the bi-focus periodic orbit $L$ contain an orbit $\Gamma$ of corank-2 homoclinic tangency. The points of intersection of the homoclinic orbit $\Gamma$ with the neighborhood $U$ belong to local stable manifold $W^s_{loc}(O)$ and local unstable manifold $W^u_{loc}(O)$, and converge to $O$ at the forward and, respectively, backward iterations of the local map $T_0$. Let $M^{+}\in W^s_{loc}(O)$ and $M^{-}\in W^u_{loc}(O)$ be two points of the orbit $\Gamma$. Since the points $M^+$ and $M^-$ belong to the the same orbit, there exists a positive integer $k_0$ such that $M^{+}=f^{k_0}(M^{-})$. Let $\Pi^{+}$ and $\Pi^{-}$ be some small neighborhoods of $M^+$ and $M^-$, respectively. The map $T_1\equiv f^{k_0}: \Pi^-\rightarrow \Pi^+$ is said to be the \textit{global map}.

If the orbit $L$ has non-leading multipliers, then we require the fulfillment of certain genericity conditions regarding the geometry of the tangency $\Gamma$. We will denote the eigensubspace of the linear part of the local map $T_0$ corresponding to the multipliers 
\begin{enumerate}

\item $\lambda_3,\ldots,\lambda_{k_1}$ (non-leading stable) by $\mathcal E^{ss}$ and $\gamma_3,\ldots,\gamma_{k_2}$ (non-leading unstable) by $\mathcal E^{uu}$;

\item $\gamma_1,\gamma_2,\lambda_1,\ldots,\lambda_{k_1}$ (leading unstable and stable) by $\mathcal E^{se}$ and $\lambda_1,\lambda_2,\gamma_1,\ldots,\gamma_{k_2}$ (leading stable and unstable) by $\mathcal E^{ue}$.

\end{enumerate}

There is an invariant $C^r$-smooth manifold $W^{ss}(O)\subset W^s(O)$ which is tangent to $\mathcal E^{ss}$. Similarly, there is an invariant $C^r$-smooth manifold $W^{uu}(O)\subset W^u(O)$, tangent to $\mathcal E^{uu}$. Recall (see e.g.  \cite{Shilnikov1}) that $W^{ss}(O)$ and $W^{uu}(O)$ are uniquely included in invariant $C^r$-smooth foliations $F^{ss}$ and $F^{uu}$ on $W^s(O)$ and $W^u(O)$, respectively. Also, there exist invariant manifolds $W^{se}(O)$ and $W^{ue}(O)$ (at least $C^1$-smooth) which are tangent to $\mathcal E^{se}$ and $\mathcal E^{ue}$, respectively. The manifolds $W^{se}(O)$ and $W^{ue}(O)$ are not unique but $\mathcal T_{M^+}\left(W^{se}(O)\right)$ and $\mathcal T_{M^-}\left(W^{ue}(O)\right)$ are defined uniquely. The genericity conditions are

\begin{doublespace}
\textbf{\hypertarget{C1.}{C1.}} $T_1\left(W^{ue}(O)\right)$ is transverse to the leaf $l^{ss}$ of $F^{ss}$ which passes through $M^+$;

\textbf{\hypertarget{C2.}{C2.}} $T_1^{-1}\left(W^{se}(O)\right)$ is transverse to the leaf $l^{uu}$ of $F^{uu}$ which passes through $M^-$.
\end{doublespace}
 
Note that since the manifolds and foliations are invariant, conditions \textbf{\hyperlink{C1.}{C1}} and \textbf{\hyperlink{C2.}{C2}} are independent of the choice of the homoclinic points $M^+$ and $M^-$. Obviously, these conditions can be fulfilled after an arbitrary $C^r$-small perturbation. A detailed discussion can be found in \cite{Turaev1}.

In the lemma below we give the formulas for the global map $T_1$ provided that the map $f$ satisfies genericity conditions \textbf{\hyperlink{C1.}{C1}} and \textbf{\hyperlink{C2.}{C2}}. Let us emphasize that these formulas define the global maps not only for the map $f$, but also for all $C^r$-small perturbations of $f$. 

\begin{lemma}\label{lemmaglobalmap}

Let a $C^r$ map $f$, where $r=2,...,\infty,\omega$, have a bi-focus periodic orbit $L$ whose stable and unstable manifolds form a corank-2 homoclinic tangency $\Gamma$ satisfying conditions \textbf{\hyperlink{C1.}{C1}},\textbf{\hyperlink{C2.}{C2}}. Then the homoclinic points $M^+$, $M^-$ and the system of coordinates in $U$ for which the local map $T_0$ is kept in the main normal form can be chosen in such a way that at all small $\varepsilon$ the Taylor expansion for the global map $T_1$ at $\left(x_1=0,x_2=0,y_1=y_1^-(\varepsilon),y_2=y_2^-(\varepsilon),u=0,\overline v=0\right)$ is as follows:
\begin{equation}\label{globalmap}
\begin{aligned} 
&\overline x_1 -x_1^+ = a_{11}x_1+a_{12}x_2+b_{11}(y_1-y_1^-)+b_{12}(y_2-y_2^-)+c_{1} u+d_{1} \overline v+\ldots, \\
&\overline x_2 -x_2^+ = a_{21}x_1+a_{22}x_2+b_{21}(y_1-y_1^-)+b_{22}(y_2-y_2^-)+c_{2} u+d_{2} \overline v+\ldots, \\
&\overline y_1 = a_{31} x_1+a_{32} x_2+\mu_{0,0}+\mu_{1,0} (y_1-y_1^-)+\mu_{1,1} (y_2-y_2^-)+c_{3} u+d_{3} \overline v+\ldots, \\
&\overline y_2=a_{41} x_1+a_{42} x_2+\nu_{0,0}+\nu_{1,0} (y_1-y_1^-)+\nu_{1,1} (y_2-y_2^-)+c_{4} u+d_{4} \overline v+\ldots, \\
&\overline u-u^+= a_{51} x_1+a_{52} x_2 + b_{51}(y_1-y_1^-)+b_{52}(y_2-y_2^-)+c_5 u+d_5\overline v+\ldots,\\
&v-v^-=a_{61} x_1+a_{62} x_2 + b_{61}(y_1-y_1^-)+b_{62}(y_2-y_2^-)+c_6 u+d_6\overline v+\ldots,\\
\end{aligned}
\end{equation}
where $\det \begin{pmatrix} a_{31}& a_{32}\\ a_{41}& a_{42}\end{pmatrix}\not=0$, $\det \begin{pmatrix} b_{11}& b_{12}\\ b_{21}& b_{22}\end{pmatrix}\not=0$, $\det d_6\not =0$ and $\mu(0)=\nu(0)=0$.

All the coefficients in \eqref{globalmap} depend on $\varepsilon$ (at least $C^{r-2}$-smoothly).
\end{lemma}

\begin{demo}

Let $M^+\in W^s_{loc}(O)$ and $M^-\in W^u_{loc}(O)$ be a pair of homoclinic points at $\varepsilon=0$. In the coordinates of Proposition \ref{proposition1} the manifolds $W^s_{loc}(O)$ and $W^u_{loc}(O)$ are straightened, so the $(y_1,y_2,v)$-coordinates of $M^+$ and $(x_1,x_2,u)$-coordinates of $M^-$ are zero. Let $M^+=M^+(x_1^+,x_2^+,0,0,u^+,0)$ and $M^-=(0,0,y_1^-,y_2^-,0,v^-)$. Since $T_1M^-=M^+$ at $\varepsilon=0$, the map $T_1(\varepsilon)$ may be written in the following form at small $\varepsilon$:
\begin{equation}\label{lemmaglobalmapeq1}
\begin{aligned} 
&\overline x_1 -x_1^+(\varepsilon) = \tilde a_{11}x_1+\tilde a_{12}x_2+\tilde b_{11}(y_1-y_1^-)+\tilde b_{12}(y_2-y_2^-)+\tilde c_{1} u+\tilde d_{1} (v-v^-)+\ldots, \\
&\overline x_2 -x_2^+(\varepsilon) = \tilde a_{21}x_1+\tilde a_{22}x_2+\tilde b_{21}(y_1-y_1^-)+\tilde b_{22}(y_2-y_2^-)+\tilde c_{2} u+\tilde d_{2} (v-v^-)+\ldots, \\
&\overline y_1 = y_1^+(\varepsilon)+\tilde a_{31} x_1+\tilde a_{32} x_2+\tilde b_{31}(y_1-y_1^-)+\tilde b_{32}(y_2-y_2^-)+\tilde c_{3} u+\tilde d_{3} (v-v^-)+\ldots, \\
&\overline y_2= y_2^+(\varepsilon)+\tilde a_{41} x_1+\tilde a_{42} x_2+\tilde b_{41}(y_1-y_1^-)+\tilde b_{42}(y_2-y_2^-)+\tilde c_{4} u+\tilde d_{4} (v-v^-)+\ldots, \\
&\overline u-u^+(\varepsilon)= \tilde a_{51} x_1+\tilde a_{52} x_2 + \tilde b_{51}(y_1-y_1^-)+\tilde b_{52}(y_2-y_2^-)+\tilde c_5 u+\tilde d_5 (v-v^-)+\ldots,\\
&\overline v=v^+(\varepsilon)+\tilde a_{61} x_1+\tilde a_{62} x_2 + \tilde b_{61}(y_1-y_1^-)+\tilde b_{62}(y_2-y_2^-)+\tilde c_6 u+\tilde d_6 (v-v^-)+\ldots,\\
\end{aligned}
\end{equation}
where the dots denote nonlinear terms, all the coefficients depend on $\varepsilon$ and $y_1^+(0)=0$, $y_2^+(0)=0$, $v^+(0)=0$. Further, for brevity, in the formulas we will omit the dependence on $\varepsilon$.

By virtue of the condition  \textbf{\hyperlink{C2.}{C2}}, the manifold $T_1^{-1}\left(W^{se}(O)\right)$ is transverse to $l^{uu}$ at the point $M^-$. Let us note that the leaf $l^{uu}$ is given by the equations $\{x_1=0,x_2=0,u=0,y_1=y_1^-,y_2=y_2^-\}$, and the tangent space $\mathcal T_{M^+}\left(W^{se}(O)\right)$ is given by $\overline v=0$. Therefore, plugging these equalities in the last line in \eqref{lemmaglobalmapeq1}, we get that the transversality condition  \textbf{\hyperlink{C2.}{C2}} is written as $\det d_6\not =0$. It means that one can resolve the last equation in \eqref{lemmaglobalmapeq1} with respect to $(v-v^-)$ and plug the result into other equations in \eqref{lemmaglobalmapeq1}. As a result, one obtains
\begin{equation}\label{lemmaglobalmapeq2}
\begin{aligned} 
&\overline x_1 -x_1^+ = a_{11}x_1+a_{12}x_2+b_{11}(y_1-y_1^-)+b_{12}(y_2-y_2^-)+c_{1} u+d_{1} \overline v+\ldots, \\
&\overline x_2 -x_2^+ = a_{21}x_1+a_{22}x_2+b_{21}(y_1-y_1^-)+b_{22}(y_2-y_2^-)+c_{2} u+d_{2} \overline v+\ldots, \\
&\overline y_1 = y_1^++a_{31} x_1+a_{32} x_2+b_{31}(y_1-y_1^-)+b_{32}(y_2-y_2^-)+c_{3} u+d_{3} \overline v+\ldots, \\
&\overline y_2= y_2^++a_{41} x_1+a_{42} x_2+b_{41}(y_1-y_1^-)+b_{42}(y_2-y_2^-)+c_{4} u+d_{4} \overline v+\ldots, \\
&\overline u-u^+= a_{51} x_1+a_{52} x_2 + b_{51}(y_1-y_1^-)+b_{52}(y_2-y_2^-)+c_5 u+d_5\overline v+\ldots,\\
&v-v^-=a_{61} x_1+a_{62} x_2 + b_{61}(y_1-y_1^-)+b_{62}(y_2-y_2^-)+c_6 u+d_6\overline v+\ldots.
\end{aligned}
\end{equation}

Since the tangency $\Gamma$ is of corank 2, then the manifolds $T_1\left(W^u_{loc}(O)\right)$ and $W^s_{loc}(O)$ have a common tangent plane at the point $M^+$ (at $\varepsilon=0$). Note that $W^u_{loc}(O)$ and $W^s_{loc}(O)$ are, respectively, given by the equations $\{x_1=0, x_2=0, u=0\}$ and $\{\overline y_1=0, \overline y_2=0, \overline v\}$. Then \eqref{lemmaglobalmapeq2} implies that the intersection of the tangent spaces to $T_1\left(W^u_{loc}(O)\right)$ and $W^s_{loc}(O)$ at the point $M^+$ is two-dimensional if and only if the following system has a two-parameter family of solutions at $\varepsilon=0$:
\begin{equation}\label{lemmaglobalmapeq4}
\begin{aligned} 
&b_{31}(y_1-y_1^-)+b_{32}(y_2-y_2^-)=0, \\
&b_{41}(y_1-y_1^-)+b_{42}(y_2-y_2^-)=0. \\
\end{aligned}
\end{equation}
This condition is fulfilled when 
\begin{equation}\label{lemmaglobalmapeq5}
b_{31}=b_{32}=b_{41}=b_{42}=0.
\end{equation}

By the condition \textbf{\hyperlink{C1.}{C1}}, the manifold $T_1\left(W^{ue}(O)\right)$ is transverse to $l^{ss}$ at the point $M^+$. The leaf $l^{ss}$ is given by the equations $\{\overline x_1=x_1^+, \overline x_2=x_2^+, \overline y_1=0, \overline y_2=0,\overline v=0\}$ and the tangent space $\mathcal T_{M^-}\left(W^{ue}(O)\right)$ is given by $\{u=0\}$. This condition can be written as
\begin{equation}\label{lemmaglobalmapeq3}
\begin{aligned} 
\det\begin{pmatrix}a_{11}&a_{12}&b_{11}&b_{12}\\a_{21}&a_{22}&b_{21}&b_{22}\\a_{31}&a_{32}&b_{31}&b_{32}\\a_{41}&a_{42}&b_{41}&b_{42}\end{pmatrix}\not=0.
\end{aligned}
\end{equation}

Combining \eqref{lemmaglobalmapeq5} and \eqref{lemmaglobalmapeq3} we get that 
\begin{equation*}
\begin{aligned} 
\det \begin{pmatrix} a_{31}& a_{32}\\ a_{41}& a_{42}\end{pmatrix}\not=0 \;\;\; \text{and} \;\;\; \det \begin{pmatrix} b_{11}& b_{12}\\ b_{21}& b_{22}\end{pmatrix}\not=0.
\end{aligned}
\end{equation*}

Finally, denoting $y_1^+=\mu^0_0$, $y_2^+=\nu^0_0$ and $b_{31}=\mu^1_0$, $b_{32}=\mu^1_1$, $b_{41}=\nu^1_0$, $b_{42}=\nu^1_1$ we bring the system \eqref{lemmaglobalmapeq2} to the form \eqref{globalmap}. 
\end{demo}

\subsection{Index of a corank-2 tangency}\label{subsection3.3}

As in the previous section, assume that the stable and unstable manifolds of the bi-focus periodic orbit $L$ contain an orbit $\Gamma$ of corank-2 homoclinic tangency. In this case, the Definition \ref{definition2} of order of tangency can be reformulated making use of a global map. Specifically, the tangency $\Gamma$ is of order $n<r$, $n\in\mathbb N$, if in the formula \eqref{globalmap} at $\varepsilon=0$ and at the point $\left(x_1=0,x_2=0,y_1=y^-_1,y_2=y^-_2,u=0,\overline v=0\right)$ one has
\begin{enumerate}

\item $\frac{\partial^{i+j} (\overline y_1,\overline y_2)}{\partial y^i_1\partial y^j_2}=0$ for all $i,j\in\{0\}\cup\mathbb N$ such that $(i+j)\in\{1,...,n\}$;

\item $\frac{\partial^{n+1} (\overline y_1,\overline y_2)}{\partial y^i_1\partial y^j_2}\not=0$ for at least one pair of $i,j\in\{0\}\cup\mathbb N$ such that $i+j=n+1$.

\end{enumerate}
If $r$ is finite and the first condition holds for $n=r$, or $r=\infty$ and the first condition holds for all $n\in\mathbb N$, then the tangency $\Gamma$ is $C^r$-flat. 

Further, we introduce \textit{suborder} of tangency $\Gamma$. Let us emphasize that, unlike the order of tangency, the suborder is not invariant under smooth changes of coordinates. Therefore, its definition is purely technical and will be used in the proofs of Lemmas \ref{lemmahsotangency} and \ref{lemmahotangency}. Again, in the following definition we employ the formula \eqref{globalmap}, where we fix $\varepsilon=0$ and the point $\left(x_1=0,x_2=0,y_1=y^-_1,y_2=y^-_2,u=0,\overline v=0\right)$.

\begin{definition}\label{definition4}
Let the tangency $\Gamma$ be order $n$. Then we say that $\Gamma$ is of suborder

\begin{enumerate}

\item $m=0$ if $\frac{\partial^{n+1} (\overline y_1,\overline y_2)}{\partial y^{n+1}_1}\not=0$;

\item $m\in\mathbb N$ if $\frac{\partial^{n+1} (\overline y_1,\overline y_2)}{\partial y^{n+1-i}_1\partial y^{i}_2}=0$ for all $0\le i<m\le n+1$, and $\frac{\partial^{n+1} (\overline y_1,\overline y_2)}{\partial y^{n+1-m}_1\partial y^m_2}\not=0$.

\end{enumerate}

\end{definition}

Thus, when the system of coordinates is fixed, we can establish 

\begin{definition}\label{definition5}
Let the tangency $\Gamma$ be order $n$ and of suborder $m$. Then we say that $\Gamma$ is of index $(n,m)$.
\end{definition}

One can introduce strict order relation for indices:  
$$(n,m)<(n',m') \iff \left[ \begin{gathered} n<n', \hfill \\ n=n' \;\;\; \text{and}\;\;\; m<m'. \hfill \\ \end{gathered} \right.$$

If the tangency $\Gamma$ is of index $(n,m)$, then for all small $\varepsilon$ the third and the fourth lines in \eqref{globalmap} (equations for $\overline y_1$ and $\overline y_2$) can be rewritten as
\begin{equation}\label{globalmap1}
\begin{aligned} 
&\overline y_1 = a_{31} x_1+a_{32} x_2+\sum\limits_{j=0}^n\sum\limits_{i=0}^{j} \mu_{j,i}(y_1-y_1^-)^{j-i} (y_2-y_2^-)^{i}+\\
&+\sum\limits_{i=0}^{m-1}\mu_{n+1,i} (y_1-y_1^-)^{n+1-i}(y_2-y_2^-)^{i}+\sum\limits_{i=m}^{n+1}A_i (y_1-y_1^-)^{n+1-i}(y_2-y_2^-)^{i}+c_{3} u+d_{3} \overline v+\ldots, \\
&\overline y_2=a_{41} x_1+a_{42} x_2+\sum\limits_{j=0}^n\sum\limits_{i=0}^{j} \nu_{j,i} (y_1-y_1^-)^{j-i} (y_2-y_2^-)^{i}+\\
&+\sum\limits_{i=0}^{m-1}\nu_{n+1,i} (y_1-y_1^-)^{n+1-i}(y_2-y_2^-)^{i}+\sum\limits_{i=m}^{n+1}B_i (y_1-y_1^-)^{n+1-i}(y_2-y_2^-)^{i}+c_{4} u+d_{4} \overline v+\ldots, \\
\end{aligned}
\end{equation}
where all the coefficients depend continuously in the $C^r$-topology, $\mu_{j,i}$ and $\nu_{j,i}$ are smooth functionals which are small in absolute value, and $|A_m|+|B_m|>0$. 

If the map $f$ is included in an arbitrary finite-parameter family of $C^r$ maps, then similarly to generical and independent splitting of arbitrary corank-2 tangency of order $n$ (Definition \ref{definition3}), we can define \textit{generical} and \textit{independent subsplitting} of the tangency $\Gamma$ of index $(n,m)$. The only difference is that the functionals $\mu_{j,i}$ and $\nu_{j,i}$ need to be taken not from \eqref{mapFperturbations}, but from \eqref{globalmap1}. Absolutely similar to the proof of Lemma \ref{lemmasplittingofcorank2tangency}, one can prove

\begin{lemma}\label{lemmasplittingoftangency}
Let a $C^r$ map $f$, where $r=2,...,\infty,\omega$, have bi-focus periodic orbits $L_1,\ldots,L_q$ whose stable and unstable manifolds form corank-2 homoclinic tangencies $\Gamma_1,\ldots,\Gamma_h$ satisfying conditions \textbf{\hyperlink{C1.}{C1}},\textbf{\hyperlink{C2.}{C2}}. Then there exists an analytic finite-parameter family of real analytic maps $\mathcal F_{\varepsilon}$ such that in the family $\mathcal F_{\varepsilon} \circ f$ the tangencies $\Gamma_1,...,\Gamma_h$ subsplit generically and independently.
\end{lemma}

\subsection{Algorithm for obtaining high order tangencies of corank 2}\label{subsection3.4}

Note that all the constructions in this section are carried out in the same coordinate system (for which the local map $T_0^k$ is in the main normal form and the global map $T_1$ satisfies \eqref{globalmap}), therefore the suborders and indices of tangencies are determined correctly. Also, let us emphasize that all the lemmas in this section hold for smooth and real-analytic cases. However, for simplicity and brevity of presentation, we formulate them only for the real-analytic case.

\begin{lemma}\label{lemmahsotangency}
Let a real-analytic map $f_0$ have a bi-focus periodic orbit $L_{f_0}$ whose stable and unstable manifolds $W^s(L_{f_0})$ and $W^u(L_{f_0})$ form two corank-2 homoclinic tangencies $\Gamma_1$, $\Gamma_2$. Let these tangencies be of index $(n,m)$ and satisfy conditions \textbf{\hyperlink{C1.}{C1}},\textbf{\hyperlink{C2.}{C2}}. Then there exists an analytic finite-parameter family of real analytic maps $f_{\varepsilon}$ and a sequence $\varepsilon_i\xrightarrow[i\rightarrow+\infty]{}0$ such that $W^s(L_{f_{\varepsilon_i}})$ and $W^u(L_{f_{\varepsilon_i}})$ form a corank-2 homoclinic tangency $\Gamma$ of index $(n,m+1)$ or $(n+1,0)$.
\end{lemma}

\begin{demo}

Let $O_{f_{0}}$ be an arbitrary point from the orbit $L_{f_{0}}$. Let $(x_1,x_2,y_1,y_2,u,v)$ be the coordinates near $O_{f_{0}}$ such that the local map $T_0$ is in the main normal form. Then, according to Proposition \ref{proposition1}, formulae \eqref{eqlocal3} hold for the map $T_0^k$. Let $T_1:\Pi_1^-\rightarrow\Pi_1^+$ and $\hat T_1:\Pi_2^-\rightarrow\Pi_2^+$ be global maps for $\Gamma_1$, $\Gamma_2$, respectively. According to Lemma \ref{lemmasplittingoftangency}, we can include the map $f_0$ in an analyitc finite-parameter family of real-analytic maps in which tangencies $\Gamma_1$ and $\Gamma_2$ subsplit generically and independently. The global maps $T_1$, $\hat T_1$ can be written as
\begin{equation}\label{lemmahsotangencyeq1}
T_1: \;\;\;
\begin{aligned} 
&\overline x_1 -x_1^+ = a_{11}x_1+a_{12}x_2+b_{11}(y_1-y_1^-)+b_{12}(y_2-y_2^-)+c_{1} u+d_{1} \overline v+\ldots, \\
&\overline x_2 -x_2^+ = a_{21}x_1+a_{22}x_2+b_{21}(y_1-y_1^-)+b_{22}(y_2-y_2^-)+c_{2} u+d_{2} \overline v+\ldots, \\
&\overline y_1 = a_{31} x_1+a_{32} x_2+\sum\limits_{j=0}^n\sum\limits_{i=0}^{j} \mu_{j,i}(y_1-y_1^-)^{j-i} (y_2-y_2^-)^{i}+\\
&\sum\limits_{i=0}^{m-1}\mu_{n+1,i} (y_1-y_1^-)^{n+1-i}(y_2-y_2^-)^{i}+\sum\limits_{i=m}^{n+1}A_i (y_1-y_1^-)^{n+1-i}(y_2-y_2^-)^{i}+c_{3} u+d_{3} \overline v+\ldots, \\
&\overline y_2=a_{41} x_1+a_{42} x_2+\sum\limits_{j=0}^n\sum\limits_{i=0}^{j} \nu_{j,i} (y_1-y_1^-)^{j-i} (y_2-y_2^-)^{i}+\\
&\sum\limits_{i=0}^{m-1}\nu_{n+1,i} (y_1-y_1^-)^{n+1-i}(y_2-y_2^-)^{i}+\sum\limits_{i=m}^{n+1}B_i (y_1-y_1^-)^{n+1-i}(y_2-y_2^-)^{i}+c_{4} u+d_{4} \overline v+\ldots, \\
&\overline u-u^+= a_{51} x_1+a_{52} x_2 + b_{51}(y_1-y_1^-)+b_{52}(y_2-y_2^-)+c_5 u+d_5\overline v+\ldots,\\
&v-v^-=a_{61} x_1+a_{62} x_2 + b_{61}(y_1-y_1^-)+b_{62}(y_2-y_2^-)+c_6 u+d_6\overline v+\ldots,\\
\end{aligned}
\end{equation}
and
\begin{equation}\label{lemmahsotangencyeq2}
\hat T_1: \;\;\;
\begin{aligned} 
&\overline{\overline{\overline x}}_1 -\hat x_1^+ = \hat a_{11}\overline{\overline x}_1+\hat a_{12}\overline{\overline x}_2+\hat b_{11}(\overline{\overline y}_1-\hat y_1^-)+\hat b_{12}(\overline{\overline y}_2-\hat y_2^-)+\hat c_{1} \overline{\overline u}+\hat d_{1} \overline{\overline{\overline v}}+\ldots, \\
&\overline{\overline{\overline x}}_2 -\hat x_2^+ = \hat a_{21}\overline{\overline x}_1+\hat a_{22}\overline{\overline x}_2+\hat b_{21}(\overline{\overline y}_1-\hat y_1^-)+\hat b_{22}(\overline{\overline y}_2-\hat y_2^-)+\hat c_{2} \overline{\overline u}+\hat d_{2} \overline{\overline{\overline v}}+\ldots, \\
&\overline{\overline{\overline y}}_1 = \hat a_{31} \overline{\overline x}_1+\hat a_{32} \overline{\overline x}_2+\sum\limits_{j=0}^n\sum\limits_{i=0}^{j} p_{j,i}(\overline{\overline y}_1-\hat y_1^-)^{j-i} (\overline{\overline y}_2-\hat y_2^-)^{i}+\\
&\sum\limits_{i=0}^{m-1}p_{n+1,i} (\overline{\overline y}_1-\hat y_1^-)^{n+1-i}(\overline{\overline y}_2-\hat y_2^-)^{i}+\sum\limits_{i=m}^{n+1}D_i (\overline{\overline y}_1-\hat y_1^-)^{n+1-i}(\overline{\overline y}_2-\hat y_2^-)^{i}+\hat c_{3} \overline{\overline u}+\hat d_{3} \overline{\overline{\overline v}}+\ldots, \\
&\overline{\overline{\overline y}}_2=\hat a_{41} \overline{\overline x}_1+\hat a_{42} \overline{\overline x}_2+\sum\limits_{j=0}^n\sum\limits_{i=0}^{j} q_{j,i} (\overline{\overline y}_1-\hat y_1^-)^{j-i} (\overline{\overline y}_2-\hat y_2^-)^{i}+\\
&\sum\limits_{i=0}^{m-1}q_{n+1,i} (\overline{\overline y}_1-\hat y_1^-)^{n+1-i}(\overline{\overline y}_2-\hat y_2^-)^{i}+\sum\limits_{i=m}^{n+1}E_i (\overline{\overline y}_1-\hat y_1^-)^{n+1-i}(\overline{\overline y}_2-\hat y_2^-)^{i}+\hat c_{4} \overline{\overline u}+\hat d_{4} \overline{\overline{\overline v}}+\ldots, \\
&\overline{\overline{\overline u}}-\hat u^+= \hat a_{51} \overline{\overline x}_1+\hat a_{52} \overline{\overline x}_2 + \hat b_{51}(\overline{\overline y}_1-\hat y_1^-)+\hat b_{52}(\overline{\overline y}_2-\hat y_2^-)+\hat c_5 \overline{\overline u}+\hat d_5\overline{\overline{\overline v}}+\ldots,\\
&\overline{\overline v}-\hat v^-=\hat a_{61} \overline{\overline x}_1+\hat a_{62} \overline{\overline x}_2 + \hat b_{61}(\overline{\overline y}_1-\hat y_1^-)+\hat b_{62}(\overline{\overline y}_2-\hat y_2^-)+\hat c_6 \overline{\overline u}+\hat d_6\overline{\overline{\overline v}}+\ldots.\\
\end{aligned}
\end{equation}
Since the tangency $\Gamma_2$ has index $(n,m)$, then, by definition, at least one of the coefficients $D_m$, $E_m$ does not equal to zero. We will assume that $D_m\not=0$.

It follows from \eqref{eqlocal3} that the equation for the local map $T_0^k$ has the following form (in formulas for $y$-coordinates we expressed $y_{k1}$, $y_{k2}$ in terms of $y_{01}$, $y_{02}$ in the linear part)
\begin{equation}\label{lemmahsotangencyeq2.5}
T_0^k: \;\;\;
\begin{aligned} 
&\overline{\overline x}_1 = \lambda^k\cdot\left(\overline x_1\cos(k\varphi)-\overline x_2\sin(k\varphi)\right)+\hat\lambda^k\cdot \xi^1_k(\overline x_1,\overline x_2,\overline{\overline y}_1,\overline{\overline y}_2,\overline u,\overline{\overline v},\mu,\nu),\\
&\overline{\overline x}_2 = \lambda^k\cdot\left(\overline x_1\sin(k\varphi)+\overline x_2\cos(k\varphi)\right)+\hat\lambda^k\cdot \xi^2_k(\overline x_1,\overline x_2,\overline{\overline y}_1,\overline{\overline y}_2,\overline u,\overline{\overline v},\mu,\nu),\\
&\overline{\overline y}_1 = \gamma^k \cdot\left(\overline y_1\cos(k\psi)-\overline y_2\sin(k\psi)\right)+\gamma^k\cdot\hat\gamma^{-k}\cdot\tilde\eta^1_k(\overline x_1,\overline x_2,\overline{\overline y}_1,\overline{\overline y}_2,\overline u,\overline{\overline v},\mu,\nu), \\
&\overline{\overline y}_2 = \gamma^k \cdot\left(\overline y_1\sin(k\psi)+\overline y_2\cos(k\psi)\right)+\gamma^k\cdot\hat\gamma^{-k}\cdot\tilde\eta^2_k(\overline x_1,\overline x_2,\overline{\overline y}_1,\overline{\overline y}_2,\overline u,\overline{\overline v},\mu,\nu), \\
&\overline{\overline u}= \hat\lambda^k\cdot \hat\xi_k(\overline x_1,\overline x_2,\overline{\overline y}_1,\overline{\overline y}_2,\overline u,\overline{\overline v},\mu,\nu),\\
&\overline v=\hat\gamma^{-k}\cdot \hat\eta_k(\overline x_1,\overline x_2,\overline{\overline y}_1,\overline{\overline y}_2,\overline u,\overline{\overline v},\mu,\nu).\\
\end{aligned}
\end{equation}
By adding arbitrarily small analytic perturbation which is localized near $L$ and does not affect tangencies $\Gamma_1$, $\Gamma_2$, we make $\varphi$ and $\psi$ rationally independent\footnote{The numbers $\varphi$, $\psi$ and 1 are said to be rationally independent if $k_1\varphi+k_2\psi$ is not an integer for any set of integers $k_1,k_2$, except $k_1=k_2=0$.} (see \cite{arbitraryorderconservative},\cite{Mints1}).

Next, we will consider a global map $T_1^{new}:\Pi^-_1\rightarrow \Pi^+_2$ such that $T^{new}_1=\hat T_1\circ T_0^k\circ T_1$ (see Fig. \ref{Figure4}) and show that there exists a sequence $\left(k_{\alpha}\right)_{\alpha\in\mathbb N}$ with the following property: $k_{\alpha}\xrightarrow[\alpha\rightarrow+\infty]{}+\infty$ and for each $\alpha$ the parameters $\mu$, $\nu$, $p$ and $q$ can be chosen in such a way that the map $T^{new}_1$ corresponds to homoclinic tangency of corank 2 of index not lower than $(n,m+1)$ or $(n+1,0)$. In addition, the values of these parameters tend to zero as $\alpha\rightarrow+\infty$. For this purpose, it is enough to consider only the third and the fourth lines in \eqref{lemmahsotangencyeq2}. Also, to simplify computations, we put $x_1,x_2,u,\overline{\overline{\overline v}}$ equal to zero, so we consider the restriction $T^{new}_1|_{(y_1,y_2)}$.

\begin{figure}[h]
\center{\includegraphics[width=0.6\linewidth]{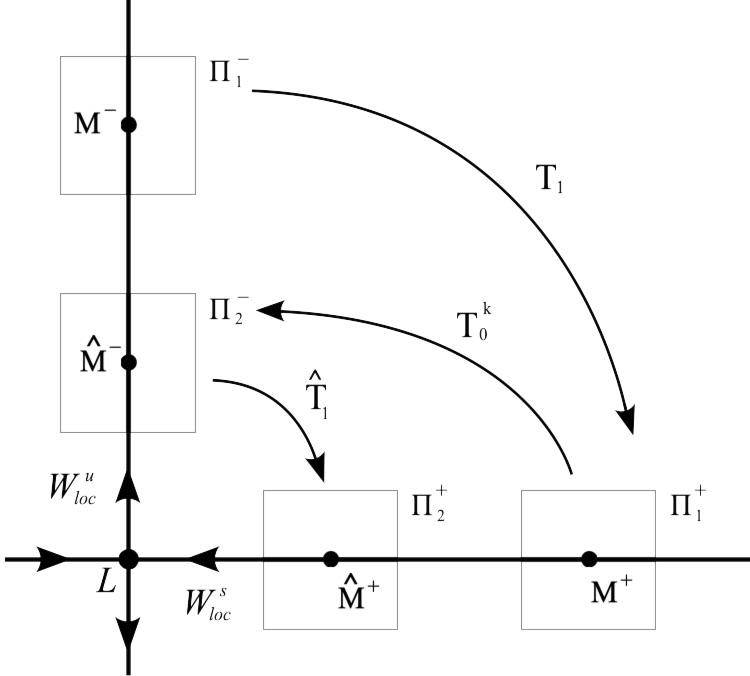}}
\caption{The global map $T_1^{new}:\Pi^-_1\rightarrow \Pi^+_2$ is given by the composition: $T^{new}_1=\hat T_1\circ T_0^k\circ T_1$, where $T_0^k$ is the local map near the periodic point $O_{f_{\varepsilon}}$ and $T_1:\Pi_1^-\rightarrow\Pi_1^+$, $\hat T_1:\Pi_2^-\rightarrow\Pi_2^+$ are the global maps for the tangencies $\Gamma_1$, $\Gamma_2$, respectively.}
\label{Figure4}
\end{figure}

Let us consider a composition $T_0^k \circ T_1|_{(y_1,y_2)}$. Combining formulas \eqref{lemmahsotangencyeq1} and \eqref{lemmahsotangencyeq2.5}, we get 
\begin{equation}\label{lemmahsotangencyeq3}
\begin{aligned} 
&\overline{\overline x}_1 = \lambda^k\cdot\left(x_1^++\tilde b_{11}Y_1+\tilde b_{12}Y_2+O(\hat\gamma^{-k})+\ldots\right)+O(\hat\lambda^k), \\
&\overline{\overline x}_2 = \lambda^k \cdot\left(x_2^++\tilde b_{21}Y_1+\tilde b_{22}Y_2+O(\hat\gamma^{-k})+\ldots\right)+O(\hat\lambda^k), \\
&\overline{\overline y}_1 = \gamma^k \cdot\left(\sum\limits_{j=0}^n\sum\limits_{i=0}^{j} \tilde\mu_{j,i}Y_1^{j-i} Y_2^{i}+\sum\limits_{i=0}^{m-1}\tilde\mu_{n+1,i} Y_1^{n+1-i}Y_2^{i}+\sum\limits_{i=m}^{n+1}\tilde A_i Y_1^{n+1-i}Y_2^{i}+O(\hat\gamma^{-k})+\ldots\right)+O(\gamma^k\cdot\hat\gamma^{-k}), \\
&\overline{\overline y}_2=\gamma^k\cdot\left(\sum\limits_{j=0}^n\sum\limits_{i=0}^{j} \tilde\nu_{j,i} Y_1^{j-i} Y_2^{i}+\sum\limits_{i=0}^{m-1}\tilde\nu_{n+1,i} Y_1^{n+1-i}Y_2^{i}+\sum\limits_{i=m}^{n+1}\tilde B_i Y_1^{n+1-i}Y_2^{i}+O(\hat\gamma^{-k})+\ldots\right)+O(\gamma^k\cdot\hat\gamma^{-k}), \\
&\overline{\overline u}= O(\hat\lambda^k),\\
&v-v^-= b_{61}Y_1+ b_{62}Y_2+ O(\hat\gamma^{-k})+\ldots,\\
\end{aligned}
\end{equation}
where we denote $Y_1=y_1-y_1^-$, $Y_2=y_2-y_2^-$ and where 
\begin{equation}\label{lemmahsotangencyeq4}
\begin{aligned} 
&\tilde b_{1i}=b_{1i}\cos(k\varphi)-b_{2i}\sin(k\varphi), \;\;\;\;\;\;\;\;\;\;\;\;\;\;\;\;\;\;\;\;\;\;\;\; \tilde b_{2i}=b_{1i}\sin(k\varphi)+b_{2i}\cos(k\varphi),\\
&\tilde\mu_{j,i}=\mu_{j,i}\cos(k\psi)-\nu_{j,i}\sin(k\psi), \;\;\;\;\;\;\;\;\;\;\;\;\;\;\;\;\;\;\;\;\; \tilde\nu_{j,i}=\mu_{j,i}\sin(k\psi)+\nu_{j,i}\cos(k\psi),\\ 
&\tilde A_m=A_m\cos(k\psi)-B_m\sin(k\psi), \;\;\;\;\;\;\;\;\;\;\;\;\;\;\;\;\;\;\;\; \tilde B_m=A_m\sin(k\psi)+B_m\cos(k\psi).\\
\end{aligned}
\end{equation}

Let us make the following change of parameters:
\begin{equation}\label{lemmahsotangencyeq5.5}
\begin{aligned} 
&\overline\mu_{j,i}=\tilde\mu_{j,i}+O(\hat\gamma^{-k}) \;\;\;\;\;\;\;\;\;\;\;\; \text{and} \;\;\; \overline\nu_{j,i}=\tilde\nu_{j,i}+O(\hat\gamma^{-k}) \;\;\;\;\;\;\;\;\;\;\;\;\; \text{for} \;\;\; 0\le j\le n, \; 0\le i\le j,\\
&\overline\mu_{n+1,i}=\tilde\mu_{n+1,i}+O(\hat\gamma^{-k}) \;\;\;\; \text{and} \;\;\; \overline\nu_{n+1,i}=\tilde\nu_{n+1,i}+O(\hat\gamma^{-k}) \;\;\;\;\; \text{for} \;\;\; 0\le i\le m-1,
\end{aligned}
\end{equation}
to get rid of the corresponding small terms (of order $O(\gamma^k\cdot\gamma^{-k})$) in the third and fourth lines.

Further, we will choose parameters $\overline\mu,\overline\nu,p,q$ to obtain the sought tangency. Let us put 
\begin{equation}\label{lemmahsotangencyeq5}
\begin{aligned} 
&\overline\mu_{0,0}=\gamma^{-k}\cdot\hat y_1^{-}, \;\;\;\;\;\;\;\;\;\;\; \overline\nu_{0,0}=\gamma^{-k}\cdot\hat y_2^{-}, \;\;\;\;\;\;\;\;\;\;\; \overline\mu_{1,1}=0, \;\;\;\;\;\;\;\;\;\;\;\overline\nu_{1,0}=0.\\
\end{aligned}
\end{equation}
In the equation \eqref{lemmahsotangencyeq2} for the map $\hat T_1$, let us put
\begin{equation}\label{lemmahsotangencyeq6}
\begin{aligned} 
&p_{j,i}=0 \;\;\;\;\;\;\; \text{and} \;\;\; q_{j,i}=0 \;\;\;\;\;\;\; \text{for} \;\;\; 2\le j\le n, \; 0\le i\le j,\\
&p_{n+1,i}=0 \;\;\; \text{and} \;\;\; q_{n+1,i}=0 \;\;\; \text{for} \;\;\; 0\le i\le m-1.
\end{aligned}
\end{equation}
Consider the composition $T^{new}_1|_{(y_1,y_2)}=\tilde T_1\circ T_0^{k}\circ T_1|_{(y_1,y_2)}$ taking into account \eqref{lemmahsotangencyeq5.5}, \eqref{lemmahsotangencyeq5} and \eqref{lemmahsotangencyeq6}. As we have already mentioned, we will write just the third and the fourth lines (for $\overline{\overline{\overline y}}_1$ and $\overline{\overline{\overline y}}_2$).
\begin{equation}\label{lemmahsotangencyeq7}
\begin{aligned} 
&\overline{\overline{\overline y}}_1=\lambda^k\cdot\hat a_{31}\left(x_1^++\tilde b_{11}Y_1+\tilde b_{12}Y_2+\ldots\right)+\lambda^k\cdot\hat a_{32}\left(x_2^++\tilde b_{21}Y_1+\tilde b_{22}Y_2+\ldots\right)+O(\hat\lambda^k)+\\
&+p_{0,0}+\gamma^k\cdot p_{1,0}\left(\overline\mu_{1,0}Y_1+\sum\limits_{j=2}^n\sum\limits_{i=0}^{j} \tilde\mu_{j,i}Y_1^{j-i} Y_2^{i}+\sum\limits_{i=0}^{m-1}\tilde\mu_{n+1,i} Y_1^{n+1-i}Y_2^{i}+\sum\limits_{i=m}^{n+1}\tilde A_i Y_1^{n+1-i}Y_2^{i}+\ldots\right)+\\
&+\gamma^k\cdot p_{1,1}\left(\overline\nu_{1,1}Y_2+\sum\limits_{j=2}^n\sum\limits_{i=0}^{j} \tilde \nu_{j,i} Y_1^{j-i} Y_2^{i}+\sum\limits_{i=0}^{m-1}\tilde\nu_{n+1,i} Y_1^{n+1-i}Y_2^{i}+\sum\limits_{i=m}^{n+1}\tilde B_i Y_1^{n+1-i}Y_2^{i}+\ldots\right)+\\
&+\gamma^{k(n+1)}\cdot\sum\limits_{h=m}^{n+1}D_h \left(\overline\mu_{1,0}Y_1+\sum\limits_{j=2}^n\sum\limits_{i=0}^{j} \tilde\mu_{j,i}Y_1^{j-i} Y_2^{i}+\sum\limits_{i=0}^{m-1}\tilde\mu_{n+1,i} Y_1^{n+1-i}Y_2^{i}+\sum\limits_{i=m}^{n+1}\tilde A_i Y_1^{n+1-i}Y_2^{i}+\ldots\right)^{n+1-h}\times\\
&\times \left(\overline\nu_{1,1}Y_2+\sum\limits_{j=2}^n\sum\limits_{i=0}^{j} \tilde \nu_{j,i} Y_1^{j-i} Y_2^{i}+\sum\limits_{i=0}^{m-1}\tilde\nu_{n+1,i} Y_1^{n+1-i}Y_2^{i}+\sum\limits_{i=m}^{n+1}\tilde B_i Y_1^{n+1-i}Y_2^{i}+\ldots\right)^{h}+O(\gamma^k\cdot\hat\gamma^{-k})+\ldots,\\
&\overline{\overline{\overline y}}_2=\lambda^k\cdot\hat a_{41}\left(x_1^++\tilde b_{11}Y_1+\tilde b_{12}Y_2+\ldots\right)+\lambda^k\cdot\hat a_{42}\left(x_2^++\tilde b_{21}Y_1+\tilde b_{22}Y_2+\ldots\right)+O(\hat\lambda^k)+\\
&+q_{0,0}+\gamma^k\cdot q_{1,0}\left(\overline\mu_{1,0}Y_1+\sum\limits_{j=2}^n\sum\limits_{i=0}^{j} \tilde\mu_{j,i}Y_1^{j-i} Y_2^{i}+\sum\limits_{i=0}^{m-1}\tilde\mu_{n+1,i} Y_1^{n+1-i}Y_2^{i}+\sum\limits_{i=m}^{n+1}\tilde A_i Y_1^{n+1-i}Y_2^{i}+\ldots\right)+\\
&+\gamma^k\cdot q_{1,1}\left(\overline\nu_{1,1}Y_2+\sum\limits_{j=2}^n\sum\limits_{i=0}^{j} \tilde \nu_{j,i} Y_1^{j-i} Y_2^{i}+\sum\limits_{i=0}^{m-1}\tilde\nu_{n+1,i} Y_1^{n+1-i}Y_2^{i}+\sum\limits_{i=m}^{n+1}\tilde B_i Y_1^{n+1-i}Y_2^{i}+\ldots\right)+\\
&+\gamma^{k(n+1)}\cdot\sum\limits_{h=m}^{n+1}E_h \left(\overline\mu_{1,0}Y_1+\sum\limits_{j=2}^n\sum\limits_{i=0}^{j} \tilde\mu_{j,i}Y_1^{j-i} Y_2^{i}+\sum\limits_{i=0}^{m-1}\tilde\mu_{n+1,i} Y_1^{n+1-i}Y_2^{i}+\sum\limits_{i=m}^{n+1}\tilde A_i Y_1^{n+1-i}Y_2^{i}+\ldots\right)^{n+1-h}\times\\
&\times \left(\overline\nu_{1,1}Y_2+\sum\limits_{j=2}^n\sum\limits_{i=0}^{j} \tilde \nu_{j,i} Y_1^{j-i} Y_2^{i}+\sum\limits_{i=0}^{m-1}\tilde\nu_{n+1,i} Y_1^{n+1-i}Y_2^{i}+\sum\limits_{i=m}^{n+1}\tilde B_i Y_1^{n+1-i}Y_2^{i}+\ldots\right)^{h}+O(\gamma^k\cdot\hat\gamma^{-k})+\ldots.\\
\end{aligned}
\end{equation}

Now, we have to choose remaining parameters $\overline\mu,\overline\nu,p,q$ in such a way that all $Y_1$, $Y_2$ terms are zero up to the suborder $Y_1^{n+1-m}\cdot Y_2^m$ (including). These values of the parameters are solutions of the system 
\begin{equation}\label{lemmahsotangencyeq8}
\begin{aligned} 
&\gamma^{k(n+1)}\cdot\overline\mu_{1,0}^{n+1-m}\cdot\overline\nu_{1,1}^m\cdot D_m+\gamma^k\cdot p_{1,0}\cdot\left(\tilde A_m+O(\hat\gamma^{-k})\right)+\gamma^k \cdot p_{1,1}\cdot\left(\tilde B_m+O(\hat\gamma^{-k})\right)=O(\lambda^k),\\
&\gamma^{k(n+1)}\cdot\overline\mu_{1,0}^{n+1-m}\cdot\overline\nu_{1,1}^m\cdot E_m+\gamma^k\cdot q_{1,0}\cdot\left(\tilde A_m+O(\hat\gamma^{-k})\right)+\gamma^k \cdot q_{1,1}\cdot\left(\tilde B_m+O(\hat\gamma^{-k})\right)=O(\lambda^k),\\
&p_{0,0}=O(\lambda^k), \;\;\;\;\;\;\;\;\;\;\;\;\;\;\;\;\;\;\;\;\;\;\;\;\;\;\;\;\;\;\;\;\;\;\;\;\;\;\;\;\;\;\;\;\;\;\;\;\;\;\;\;\;\;\;\;\;\;\;\;\;\;\;\;\;\;\; q_{0,0}=O(\lambda^k),\\
&\gamma^k\cdot p_{1,0}\cdot\overline\mu_{1,0}=\lambda^k\cdot S_{1}+O(\hat\lambda^k),  \;\;\;\;\;\;\;\;\;\;\;\;\;\;\;\;\;\;\;\;\;\;\;\;\;\;\;\;\;\;\;\;\;\;\;\;\;  \gamma^k\cdot q_{1,0}\cdot\overline\mu_{1,0}=\lambda^k \cdot S_{3}+O(\hat\lambda^k),\\
&\gamma^k\cdot p_{1,1}\cdot\overline\nu_{1,1}=\lambda^k\cdot S_{2}+O(\hat\lambda^k),  \;\;\;\;\;\;\;\;\;\;\;\;\;\;\;\;\;\;\;\;\;\;\;\;\;\;\;\;\;\;\;\;\;\;\;\;\;\;  \gamma^k\cdot q_{1,1}\cdot\overline\nu_{1,1}=\lambda^k \cdot S_{4}+O(\hat\lambda^k),\\
&\gamma^k\left(p_{1,0}\cdot\overline\mu_{j,i}+p_{1,1}\cdot\overline\nu_{j,i}\right)=O(\lambda^k),\;\;\;\;\;\;\;\;\;\;\;\;\;\;\;\;\;\;\;\;\;\;\;\;\;\;\;\;\;\;\;\;\;\;\;   \gamma^k\left(q_{1,0}\cdot\overline\mu_{j,i}+q_{1,1}\cdot\overline\nu_{j,i}\right)=O(\lambda^k),
\end{aligned}
\end{equation}
where $2\le j\le n, \; 0\le i\le j$ and $j=n+1,\; 0\le i\le m-1$, and where $S_1,S_2,S_3,S_4$ are constants which have the form
\begin{equation}\label{lemmahsotangencyeq9}
\begin{aligned} 
&S_1= \hat a_{31}\tilde b_{11}+\hat a_{32}\tilde b_{21}=\hat a_{31}\left(b_{11}\cos(k\varphi)-b_{21}\sin(k\varphi)\right)+\hat a_{32}\left(b_{11}\sin(k\varphi)+b_{21}\cos(k\varphi)\right),\\
&S_2= \hat a_{31}\tilde b_{12}+\hat a_{32}\tilde b_{22}=\hat a_{31}\left(b_{12}\cos(k\varphi)-b_{22}\sin(k\varphi)\right)+\hat a_{32}\left(b_{12}\sin(k\varphi)+b_{22}\cos(k\varphi)\right),\\
&S_3= \hat a_{41}\tilde b_{11}+\hat a_{42}\tilde b_{21}=\hat a_{41}\left(b_{11}\cos(k\varphi)-b_{21}\sin(k\varphi)\right)+\hat a_{42}\left(b_{11}\sin(k\varphi)+b_{21}\cos(k\varphi)\right),\\
&S_4= \hat a_{41}\tilde b_{12}+\hat a_{42}\tilde b_{22}=\hat a_{41}\left(b_{12}\cos(k\varphi)-b_{22}\sin(k\varphi)\right)+\hat a_{42}\left(b_{12}\sin(k\varphi)+b_{22}\cos(k\varphi)\right).\\
\end{aligned}
\end{equation}

Let us make rescaling of the parameters
\begin{equation}\label{lemmahsotangencyeq10}
\begin{aligned} 
&M_{1,0}=\lambda^{-\frac{k}{n+2}}\cdot\gamma^{\frac{k(n+1)}{n+2}}\cdot\overline\mu_{1,0}, \;\;\;\;\;\;\;\;\;\;\;\;\;\;\;\;\;\;\;\;\;\;\;\;\;\;\;\;\;\;\;\;\;\;\; N_{1,1}=\lambda^{-\frac{k}{n+2}}\cdot\gamma^{\frac{k(n+1)}{n+2}}\cdot\overline\nu_{1,1},\\
&M_{j,i}=\lambda^{-\frac{k}{n+2}}\cdot\gamma^{\frac{k(n+1)}{n+2}}\cdot\overline\mu_{j,i}, \;\;\;\;\;\;\;\;\;\;\;\;\;\;\;\;\;\;\;\;\;\;\;\;\;\;\;\;\;\;\;\;\;\;\;\; N_{j,i}=\lambda^{-\frac{k}{n+2}}\cdot\gamma^{\frac{k(n+1)}{n+2}}\cdot\overline\nu_{j,i},\\
&P_{0,0}=\lambda^{-k}\cdot p_{0,0}, \;\;\;\;\;\;\;\;\;\;\;\;\;\;\;\;\;\;\;\;\;\;\;\;\;\;\;\;\;\;\;\;\;\;\;\;\;\;\;\;\;\;\;\;\;\;\;\;\;\;\;\;\; Q_{0,0}=\lambda^{-k}\cdot q_{0,0},\\
&P_{1,0}=\lambda^{-\frac{k(n+1)}{n+2}}\cdot\gamma^{\frac{k}{n+2}}\cdot p_{1,0}, \;\;\;\;\;\;\;\;\;\;\;\;\;\;\;\;\;\;\;\;\;\;\;\;\;\;\;\;\;\;\;\;\;\;\;\; Q_{1,0}=\lambda^{-\frac{k(n+1)}{n+2}}\cdot\gamma^{\frac{k}{n+2}}\cdot q_{1,0},\\
&P_{1,1}=\lambda^{-\frac{k(n+1)}{n+2}}\cdot\gamma^{\frac{k}{n+2}}\cdot p_{1,1}, \;\;\;\;\;\;\;\;\;\;\;\;\;\;\;\;\;\;\;\;\;\;\;\;\;\;\;\;\;\;\;\;\;\;\;\; Q_{1,1}=\lambda^{-\frac{k(n+1)}{n+2}}\cdot\gamma^{\frac{k}{n+2}}\cdot q_{1,1}.\\
\end{aligned}
\end{equation}
After it, the system \eqref{lemmahsotangencyeq8} recasts as 
\begin{equation}\label{lemmahsotangencyeq11}
\begin{aligned} 
&M_{1,0}^{n+1-m}\cdot N_{1,1}^m\cdot D_m+P_{1,0}\cdot\left(\tilde A_m+O(\hat\gamma^{-k})\right)+P_{1,1}\cdot\left(\tilde B_m+O(\hat\gamma^{-k})\right)=O\left(\lambda^{\frac{k}{n+2}}\cdot\gamma^{-\frac{k(n+1)}{n+2}}\right),\\
&M_{1,0}^{n+1-m}\cdot N_{1,1}^m\cdot E_m+Q_{1,0}\cdot\left(\tilde A_m+O(\hat\gamma^{-k})\right)+Q_{1,1}\cdot\left(\tilde B_m+O(\hat\gamma^{-k})\right)=O\left(\lambda^{\frac{k}{n+2}}\cdot\gamma^{-\frac{k(n+1)}{n+2}}\right),\\
&P_{0,0}=O(1), \;\;\;\;\;\;\;\;\;\;\;\;\;\;\;\;\;\;\;\;\;\;\;\;\;\;\;\;\;\;\;\;\;\;\;\;\;\;\;\;\;\;\;\;\;\;\;\;\;\;\;\;\;\;\;\;\;\;\;\;\;\;\;  Q_{0,0}=O(1),\\
&P_{1,0}=\frac{S_{1}}{M_{1,0}}+o(1)_{k\to+\infty},  \;\;\;\;\;\;\;\;\;\;\;\;\;\;\;\;\;\;\;\;\;\;\;\;\;\;\;\;\;\;\;\;\;\;\;\;\;\;\;\;\;\; Q_{1,0}=\frac{S_{3}}{M_{1,0}}+o(1)_{k\to+\infty},\\
&P_{1,1}=\frac{S_2}{N_{1,1}}+o(1)_{k\to+\infty},  \;\;\;\;\;\;\;\;\;\;\;\;\;\;\;\;\;\;\;\;\;\;\;\;\;\;\;\;\;\;\;\;\;\;\;\;\;\;\;\;\;\; Q_{1,1}=\frac{S_4}{N_{1,1}}+o(1)_{k\to+\infty},\\
&P_{1,0}\cdot M_{j,i}+P_{1,1}\cdot N_{j,i}=O(1),\;\;\;\;\;\;\;\;\;\;\;\;\;\;\;\;\;\;\;\;\;\;\;\;\;\;\;\;\;\;\;\;\;\;  Q_{1,0}\cdot M_{j,i}+Q_{1,1}\cdot N_{j,i}=O(1).
\end{aligned}
\end{equation}

Next, we rewrite the system \eqref{lemmahsotangencyeq11} as
\begin{equation}\label{lemmahsotangencyeq12}
\begin{aligned} 
&M_{1,0}^{n+2-m}\cdot N_{1,1}^{m+1}\cdot D_m+S_1\cdot N_{1,1}\cdot\tilde A_m+S_2\cdot M_{1,0}\cdot\tilde B_m=o(1)_{k\to+\infty},\\
&M_{1,0}^{n+2-m}\cdot N_{1,1}^{m+1}\cdot E_m+S_3\cdot N_{1,1}\cdot\tilde A_m+S_4\cdot M_{1,0}\cdot\tilde B_m=o(1)_{k\to+\infty},\\
&P_{0,0}=O(1), \;\;\;\;\;\;\;\;\;\;\;\;\;\;\;\;\;\;\;\;\;\;\;\;\;\;\;\;\;\;\;\;\;\;\;\;\;\;\;\;\;\;\;\;\;\;\;\;\;\;\;\;\;\;\;\;\;\;\;\;\;\;\;  Q_{0,0}=O(1),\\
&P_{1,0}=\frac{S_{1}}{M_{1,0}}+o(1)_{k\to+\infty},  \;\;\;\;\;\;\;\;\;\;\;\;\;\;\;\;\;\;\;\;\;\;\;\;\;\;\;\;\;\;\;\;\;\;\;\;\;\;\;\;\;\; Q_{1,0}=\frac{S_{3}}{M_{1,0}}+o(1)_{k\to+\infty},\\
&P_{1,1}=\frac{S_2}{N_{1,1}}+o(1)_{k\to+\infty},  \;\;\;\;\;\;\;\;\;\;\;\;\;\;\;\;\;\;\;\;\;\;\;\;\;\;\;\;\;\;\;\;\;\;\;\;\;\;\;\;\;\;\; Q_{1,1}=\frac{S_4}{N_{1,1}}+o(1)_{k\to+\infty},\\
&S_1\cdot N_{1,1}\cdot M_{j,i}+S_2\cdot M_{1,0}\cdot N_{j,i}=O(1),\;\;\;\;\;\;\;\;\;\;\;\;\;\;\;\;\;\;\;\; S_3\cdot N_{1,1}\cdot M_{j,i}+S_4\cdot M_{1,0}\cdot N_{j,i}=O(1).
\end{aligned}
\end{equation}

According to Lemma \ref{lemmaglobalmap}, $\det \begin{pmatrix} b_{11}& b_{12}\\ b_{21}& b_{22}\end{pmatrix}\not=0$ and $\det \begin{pmatrix}\hat a_{31}& \hat a_{32}\\ \hat a_{41}& \hat a_{42}\end{pmatrix}\not=0$. It is directly verified that it implies 
\begin{equation}\label{lemmahsotangencyeq13}
\begin{aligned} 
&S_1S_4-S_2S_3\not=0.\\
\end{aligned}
\end{equation}
Also, taking into account these conditions, the facts that $|A_m|+|B_m|\not=0, D_m\not=0$ and $\varphi,\psi$ are rationally independent, we can choose a sequence $\left(k_{\beta}\right)_{\beta\in\mathbb N}$ such that $k_{\beta}\xrightarrow[\beta\rightarrow+\infty]{}+\infty$ and
\begin{equation}\label{lemmahsotangencyeq14}
\begin{aligned} 
&S_2E_m-S_4D_m\not=0, \;\;\;\;\;\;\;\;\;\;\;\;\;\;\;\; S_1E_m-S_3D_m\not=0,\\
&\tilde A_m>0, \tilde B_m>0 \;\;\;\;\;\;\;\;\;\;\;\;\;\;\;\;\;\;\;\;\;\; \text{for even $\beta$},\\
&\tilde A_m<0, \tilde B_m<0 \;\;\;\;\;\;\;\;\;\;\;\;\;\;\;\;\;\;\;\;\;\; \text{for odd $\beta$}.\\
\end{aligned}
\end{equation}

Now, making use of the conditions \eqref{lemmahsotangencyeq13} and \eqref{lemmahsotangencyeq14}, we can easily solve the system \eqref{lemmahsotangencyeq12} which becomes linear after the expression of the nonlinear term $M_{1,0}^{n+2-m}\cdot N_{1,1}^{m+1}$ in the first line and subsequent substitution of it to the second line. The solution has the form
\begin{equation}\label{lemmahsotangencyeq15}
\begin{aligned} 
&M_{1,0}=\sqrt[n+2]{\frac{\left(\tilde A_m\right)^{m+1}}{\left(\tilde B_m\right)^m}\cdot \tilde S_1}+o(1)_{k_{\alpha}\to+\infty}, \;\;\;\;\;\;\; N_{1,1}=\tilde S_2\cdot \sqrt[n+2]{\frac{\left(\tilde B_m\right)^{n+2-m}}{\left(\tilde A_m\right)^{n+1-m}}\cdot \tilde S_1}+o(1)_{k_{\alpha}\to+\infty},\\
&P_{0,0}=O(1), \;\;\;\;\;\;\;\;\;\;\;\;\;\;\;\;\;\;\;\;\;\;\;\;\;\;\;\;\;\;\;\;\;\;\;\;\;\;\;\;\;\;\;\;\;\;\;\;\;\;\;\;\; Q_{0,0}=O(1),\\
&P_{1,0}=\frac{S_{1}}{M_{1,0}}+o(1)_{k_{\alpha}\to+\infty},  \;\;\;\;\;\;\;\;\;\;\;\;\;\;\;\;\;\;\;\;\;\;\;\;\;\;\;\;\;\;\; Q_{1,0}=\frac{S_{3}}{M_{1,0}}+o(1)_{k_{\alpha}\to+\infty},\\
&P_{1,1}=\frac{S_2}{N_{1,1}}+o(1)_{k_{\alpha}\to+\infty},  \;\;\;\;\;\;\;\;\;\;\;\;\;\;\;\;\;\;\;\;\;\;\;\;\;\;\;\;\;\;\;\; Q_{1,1}=\frac{S_4}{N_{1,1}}+o(1)_{k_{\alpha}\to+\infty},\\
&M_{j,i}=O(1), \;\;\;\;\;\;\;\;\;\;\;\;\;\;\;\;\;\;\;\;\;\;\;\;\;\;\;\;\;\;\;\;\;\;\;\;\;\;\;\;\;\;\;\;\;\;\;\;\;\;\;\;\; N_{j,i}=O(1),
\end{aligned}
\end{equation}
where 
\begin{equation*}
\begin{aligned} 
&\tilde S_1 = \left(S_2 S_3-S_1 S_4\right)\cdot\frac{\left(S_1E_m-S_3D_m\right)^m}{\left(S_2E_m-S_4D_m\right)^{m+1}}\not=0, \;\;\;\;\;\;\;\;\;\;\;\;\;\; \tilde S_2=-\frac{S_2E_m-S_4D_m}{S_1E_m-S_3D_m}\not=0,\\
\end{aligned}
\end{equation*}

and where $\left(k_{\alpha}\right)_{\alpha\in\mathbb N}$ is a subsequence of $\left(k_{\beta}\right)_{\beta\in\mathbb N}$ corresponding to odd or even values of $\beta$. We choose it as follows. If $n+2$ is an odd number, the subsequence can be chosen arbitrarily: for both even and odd values of $\beta$ (since the system \eqref{lemmahsotangencyeq11} always has a non-zero solution in this case). If $n+2$ is an even number, we choose the subsequence $\left(k_{\alpha}\right)_{\alpha\in\mathbb N}$ depending on the sign of the constant $\tilde S_1$ in order to make the expression under the root positive (see \eqref{lemmahsotangencyeq15}). This choice of subsequence guarantees the existence of a non-zero solution of the system \eqref{lemmahsotangencyeq11}. Since the corresponding values of $M,N,P,Q$ are bounded for all large $k_{\alpha}$, the respective values of $\mu,\nu,p,q$ tend to zero as $\alpha\to+\infty$ (see \eqref{lemmahsotangencyeq10}).

Thus, the found non-zero solution proves that the map $T^{new}_1$ has homoclinic tangency of corank 2 of index not lower than $(n,m+1)$ or $(n+1,0)$. By adding arbitrarily small local analytic perturbation given by Lemma \ref{lemmasplittingoftangency}, we make this tangency of index $(n,m+1)$ or $(n+1,0)$. It completes the proof.
\end{demo}

\begin{lemma}\label{lemmahotangency}
Let a real-analytic map $f_0$ have a bi-focus periodic orbit $L_{f_0}$ whose stable and unstable manifolds $W^s(L_{f_0})$ and $W^u(L_{f_0})$ form $h=2^{n+2}$ corank-2 homoclinic tangencies $\Gamma_1,\ldots,\Gamma_h$. Let these tangencies be of index $(n,0)$ and satisfy conditions \textbf{\hyperlink{C1.}{C1}},\textbf{\hyperlink{C2.}{C2}}. Then there exists an analytic finite-parameter family of real-analytic maps $f_{\varepsilon}$ and a sequence $\varepsilon_i\xrightarrow[i\rightarrow+\infty]{}0$ such that $W^s(L_{f_{\varepsilon_i}})$ and $W^u(L_{f_{\varepsilon_i}})$ form a corank-2 homoclinic tangency $\Gamma_{1,\ldots,h}$ of index $(n+1,0)$.
\end{lemma}

\begin{demo}

To obtain a map with the sought tangency $\Gamma_{1,\ldots,h}$, we provide an algorithm which consists of $n+2$ steps. At every next step, we add a finite number of arbitrarily $C^{\omega}$-small perturbations to the map obtained at the previous step in order to construct a new map with homoclinic tangencies of higher suborder (compared to the previous step). Each perturbation that we use is given by Lemma \ref{lemmahsotangency}. Moreover, each perturbation is localized in small neighborhoods of homoclinic tangencies under consideration therefore it does not affect any given finite number of other tangencies. 

\textbf{Step 1.}
Perturbing the map $f_0$ in small neighborhoods of the tangencies $\Gamma_1$ and $\Gamma_2$, one gets a map $f_{1,1}$ with a corank-2 tangency $\Gamma_{1,2}$ of index $(n,1)$ and corank-2 tangencies $\Gamma_3,\ldots,\Gamma_h$ of index $(n,0)$. Next, perturbing the map $f_{1,1}$ in small neighborhoods of the tangencies $\Gamma_3$ and $\Gamma_4$, one gets a map $f_{1,2}$ with corank-2 tangencies $\Gamma_{1,2},\Gamma_{3,4}$ of index $(n,1)$ and corank-2 tangencies $\Gamma_5,\ldots,\Gamma_h$ of index $(n,0)$. Repeating this $2^{n+1}$ times, one gets a map $f_1\equiv f_{1,2^{n+1}}$ with $2^{n+1}$ corank-2 tangencies $\Gamma_{1,2},\ldots,\Gamma_{h-1,h}$  of index $(n,1)$.

\textbf{Step 2.}
Perturbing the map $f_1$ in small neighborhoods of the tangencies $\Gamma_{1,2}$ and $\Gamma_{3,4}$, one gets a map $f_{2,1}$ with a corank-2 tangency $\Gamma_{1,2,3,4}$ of index $(n,2)$ and corank-2 tangencies $\Gamma_{5,6},\ldots,\Gamma_{h-1,h}$ of index $(n,1)$. 
Further, perturbing the map $f_{2,1}$ in small neighborhoods of the tangencies $\Gamma_{5,6}$ and $\Gamma_{7,8}$, one gets a map $f_{2,2}$ with corank-2 tangencies $\Gamma_{1,2,3,4},\Gamma_{5,6,7,8}$ of index $(n,2)$ and corank-2 tangencies $\Gamma_9,\ldots,\Gamma_h$ of index $(n,1)$. Continuing this procedure, one gets a map $f_2\equiv f_{2,2^n}$ with $2^{n}$ corank-2 tangencies $\Gamma_{1,2,3,4},\ldots,\Gamma_{h-3,h-2,h-1,h}$ of index $(n,2)$.

Reasoning similarly, let us perform $n+2$ steps of this algorithm. For clarity, we describe the last step.

\textbf{Step n+2.}
After $(n+1)$-th step, we have a map $f_{n+1}$ with two corank-2 tangencies of index $(n,n+1)$: tangency $\Gamma_{1,\ldots,2^{n+1}}$ and tangency $\Gamma_{2^{n+1},\ldots,2^{n+2}}$. Perturbing the map $f_{n+1}$ in small neighborhoods of these tangencies, one gets a map $f_{n+2}$ with a corank-2 tangency $\Gamma_{1,\ldots,h}=\Gamma_{1,\ldots,2^{n+2}}$ of index $(n+1,0)$.

\end{demo}

\begin{lemma}\label{lemmahotangency1}
Let a real-analytic map $f_0$ have a bi-focus periodic orbit $L_{f_0}$ whose stable and unstable manifolds $W^s(L_{f_0})$ and $W^u(L_{f_0})$ form $h=2^{\frac{(N-1)(N+4)}{2}}$ corank-2 homoclinic tangencies. Let these tangencies be of index $(1,0)$ and satisfy conditions \textbf{\hyperlink{C1.}{C1}},\textbf{\hyperlink{C2.}{C2}}. Then there exists an analytic finite-parameter family of real-analytic maps $f_{\varepsilon}$ and a sequence $\varepsilon_i\xrightarrow[i\rightarrow+\infty]{}0$ such that $W^s(L_{f_{\varepsilon_i}})$ and $W^u(L_{f_{\varepsilon_i}})$ form a corank-2 homoclinic tangency of index $(N,0)$.
\end{lemma}

\begin{demo}

Note that if $N=1$, then the map $f_0$ has corank-2 homoclinic tangency of index $(1,0)$, so the lemma is proved. Let $N>1$.

To construct a map with the required properties, we specify an algorithm consisting of $N-1$ steps. When performing each next step, we give a finite number of arbitrarily $C^{\omega}$-small perturbations in order to construct a new map with homoclinic tangencies of higher order (compared to the map from the previous step). Each perturbation is provided by Lemma \ref{lemmahotangency} and is localized in small neighborhoods of homoclinic tangencies under consideration therefore it does not affect any given finite number of other tangencies. 

\textbf{Step 1.} Let us divide all given tangencies into groups of eight elements (in an arbitrary manner). Successively applying Lemma \ref{lemmahotangency} to each such group, we obtain a map $f_1$ with $h_1=\frac{h}{8}$ corank-2 homoclinic tangencies of index $(2,0)$. If $N=2$, then $h_1=1$, so the map $f_1$ has corank-2 homoclinic tangency of index $(2,0)$ and the lemma is proved. Let $N>2$. 

\textbf{Step 2.} 
Again, we arbitrarily divide obtained corank-2 tangencies of index $(2,0)$, but now into groups of sixteen elements. Applying Lemma \ref{lemmahotangency} to each such group in turn, provides a map $f_2$ with $h_2=\frac{h_1}{16}=\frac{h}{128}$ corank-2 homoclinic tangencies of index $(3,0)$. If $N=3$, then $h_2=1$, so the map $f_2$ has corank-2 homoclinic tangency of index $(3,0)$ and the lemma is proved. Let $N>3$. 

Making use of the same arguments, let us perform $N-1$ step of this procedure. Let us describe the last step.

\textbf{Step N--1.} 
After the previous step, we have a map $f_{N-2}$ with $h_{N-2}=\frac{h_{N-3}}{2^N}$ corank-2 homoclinic tangencies of index $(N-1,0)$. Dividing all these tangencies into groups of $2^{N+1}$ elements and successively applying Lemma \ref{lemmahotangency} to each such group, we obtain a map $f_{N-1}$ with $h_{N-1}=\frac{h_{N-2}}{2^{N+1}}=1$ corank-2 homoclinic tangencies of index $(N,0)$. 

It completes the proof of the lemma.
\end{demo}

\begin{demo2}

Let $f$ be a $C^r$ map $(r=2,\ldots,\infty,\omega)$ having a bi-focus periodic orbit $L_f$ whose stable and unstable manifolds form $2^{\frac{(N-1)(N+4)}{2}}$ homoclinic tangencies of corank 2. We will consider that the map $f$ is real-analytic (if $f$ is only smooth, then arbitrarily close to $f$ one can always find a map $f^{new}$ which is real-analytic and has a bi-focus periodic orbit $L_{f^{new}}$ whose stable and unstable manifolds form $2^{\frac{(N-1)(N+4)}{2}}$ homoclinic tangencies of corank 2). Making use of Lemma \ref{lemmasplittingoftangency}, near each given tangency we add small local analytic perturbation making it of index $(1,0)$. Further, applying Lemma \ref{lemmahotangency1}, arbitrarily $C^{\omega}$-close we find a map $g$ with corank-2 homoclinic tangency of order $N$ between stable and unstable manifolds of $L_g$.
\end{demo2}

\section{Universal 2-dimensional dynamics}\label{section4}

Let $f$ be a map which have a bi-focus periodic orbit $L$ whose stable and unstable manifolds form a homoclinic tangency. Let $\Pi^{+}$ and $\Pi^{-}$ be small neighborhoods of the homoclinic points $M^+\in W^s_{loc}(L)$ and $M^-\in W^u_{loc}(L)$, respectively, which are defined in Section \ref{subsection3.2}. By \eqref{eqlocal3}, the piece of $W^s_{loc}(L)$ near $M^+$ is the limit of the countable sequence of strips $\sigma_k=\Pi^+\cap T_0^{-k}\Pi^-$ on 
which the \textit{first-return maps} $T_k\equiv T_1\circ T_0^k$ are defined (see the discussion in \cite{Shilnikov1}, sections 3.8. and 3.9.). The map $T_k$ takes the strip $\sigma_k\subset\Pi^+$ back into $\Pi^+$. The same is true for all maps close to $f$.

In order to prove Theorem \ref{theorem4}, we need explicit formulas for the first-return maps $T_{k}$ near corank-2 homoclinic tangencies of high orders. To derive it, we make use of the following form of the global map $T_1$ for corank-2 tangency of order $n$ (see Section \ref{subsection3.2}).
\begin{equation}\label{lemmafirst-returnmapglobalmap}
T_1: \;\;\;
\begin{aligned} 
&\overline x_1 -x_1^+ = a_{11}x_1+a_{12}x_2+b_{11}(y_1-y_1^-)+b_{12}(y_2-y_2^-)+c_{1} u+d_{1} \overline v+\ldots, \\
&\overline x_2 -x_2^+ = a_{21}x_1+a_{22}x_2+b_{21}(y_1-y_1^-)+b_{22}(y_2-y_2^-)+c_{2} u+d_{2} \overline v+\ldots, \\
&\overline y_1 = a_{31} x_1+a_{32} x_2+\sum\limits_{j=0}^n\sum\limits_{i=0}^{j} \mu_{j,i}(y_1-y_1^-)^{j-i} (y_2-y_2^-)^{i}+\\
&+\sum\limits_{i=0}^{n+1}A_i (y_1-y_1^-)^{n+1-i}(y_2-y_2^-)^{i}+c_{3} u+d_{3} \overline v+\ldots, \\
&\overline y_2=a_{41} x_1+a_{42} x_2+\sum\limits_{j=0}^n\sum\limits_{i=0}^{j} \nu_{j,i} (y_1-y_1^-)^{j-i} (y_2-y_2^-)^{i}+\\
&+\sum\limits_{i=0}^{n+1}B_i (y_1-y_1^-)^{n+1-i}(y_2-y_2^-)^{i}+c_{4} u+d_{4} \overline v+\ldots, \\
&\overline u-u^+= a_{51} x_1+a_{52} x_2 + b_{51}(y_1-y_1^-)+b_{52}(y_2-y_2^-)+c_5 u+d_5\overline v+\ldots,\\
&v-v^-=a_{61} x_1+a_{62} x_2 + b_{61}(y_1-y_1^-)+b_{62}(y_2-y_2^-)+c_6 u+d_6\overline v+\ldots,\\
\end{aligned}
\end{equation}
where $\sum\limits_{i=0}^{n+1} |A_i|+\sum\limits_{i=0}^{n+1} |B_i|>0$.

\begin{lemma}\label{lemmafirst-returnmap}
Let a map $f_0$ have a bi-focus periodic orbit $L_{f_0}$ whose stable and unstable manifolds form a corank-2 homoclinic tangency $\Gamma$ which is of order $n$ and satisfies conditions \textbf{\hyperlink{C1.}{C1}},\textbf{\hyperlink{C2.}{C2}}. Let $f_0$ be embedded into a finite-paramter family of maps $f_{\varepsilon}$ of class $C^r$ $(r=n+1,...,\infty,\omega)$ with respect to the coordinates and the parameters $\varepsilon$ such that the tangency $\Gamma$ splits generically as $\varepsilon$ varies. Then for all sufficiently large $k\in\mathbb N$ and all sufficiently small $\varepsilon$ there exist a $C^{r}$-smooth coordinate transformation on $\sigma_{k}$: $(x,y,u,v)\rightarrow(X,Y,U,V)$ and a $C^{r-\max(2,n)}$-smooth transformation of the parameters $\varepsilon\rightarrow(M,N)$ which bring the first-return map $T_{k}$ to the following form:
\begin{equation}\label{lemmafirst-returnmapeq0}
\begin{aligned} 
&\overline Y_1=\sum\limits_{j=0}^n\sum\limits_{i=0}^{j} M_{j,i} Y_1^{j-i} Y_2^{i}+\sum\limits_{i=0}^{n+1} \tilde A_i Y_1^{n+1-i}Y_2^{i}+o(1)_{k\rightarrow+\infty},\\
&\overline Y_2=\sum\limits_{j=0}^n\sum\limits_{i=0}^{j} N_{j,i} Y_1^{j-i} Y_2^{i}+\sum\limits_{i=0}^{n+1} \tilde B_i Y_1^{n+1-i}Y_2^{i}+ o(1)_{k\rightarrow+\infty},\\
&\left(\overline X_1,\overline X_2,\overline U,V\right)=o(1)_{k\rightarrow+\infty} \\
\end{aligned}
\end{equation}
with $\tilde A_i=A_i\cos(k\psi)-B_i\sin(k\psi)$ and $\tilde B_i=A_i\sin(k\psi)+B_i\cos(k\psi)$, where $i\in\{0,\ldots,n+1\}$ and $A_i, B_i$ are given by \eqref{lemmafirst-returnmapglobalmap}.

The range of values of the new variables $(X,Y,U,V)$ and parameters $M, N$ covers a centered-at-zero ball whose radius tends to infinity as $k$ increases. The $o(1)$-terms tend to zero uniformly on any compact, along with the derivatives with respect to $(X,Y,U,V)$ up to the order $r$ and with respect to $M,N$ up to the order $\left(r-\max(2,n)\right)$.

\end{lemma}

\begin{demo}

Tangency $\Gamma$ splits generically as $\varepsilon$ varies, so the condition \eqref{corank2splitgenerically} allows us to choose $\mu$ and $\nu$ as new parameters. Further, making use of the formulas \eqref{eqlocal3} and \eqref{lemmafirst-returnmapglobalmap} for the maps $T_0^k$ and $T_1$, respectively, we get the following formula for the first-return map $T_{k}$ for all sufficiently large $k$ and all small $\varepsilon$:
\begin{equation}\label{lemmafirst-returnmapeq1}
\begin{aligned} 
&\overline x_{01}-x_1^+=\lambda^k \alpha_{11}(k\varphi)x_{01}+\hat\lambda^ka_{11}\xi^1_k+\lambda^k \alpha_{12}(k\varphi)x_{02}+\hat\lambda^ka_{12}\xi^2_k+b_{11} (y_{k1}-y_1^-)+b_{12}(y_{k2}-y_2^-)+\\
&+c_1\hat\lambda^k\hat\xi_k+d_1\hat\gamma^{-k}\hat\eta_k+\ldots,\\
&\overline x_{02}-x_2^+=\lambda^k \alpha_{21}(k\varphi)x_{01}+\hat\lambda^ka_{21}\xi^1_k+\lambda^k \alpha_{22}(k\varphi)x_{02}+\hat\lambda^ka_{22}\xi^2_k+b_{21}(y_{k1}-y_1^-)+b_{22}(y_{k2}-y_2^-)+\\
&+c_2\hat\lambda^k\hat\xi_k+d_2\hat\gamma^{-k}\hat\eta_k+\ldots,\\
&\gamma^{-k}\cdot\left(\overline y_{k1}\cos(k\psi)+\overline y_{k2}\sin(k\psi)\right)+\hat\gamma^{-k}\eta^1_k=\lambda^k \alpha_{31}(k\varphi)x_{01}+\hat\lambda^ka_{31}\xi^1_k+\lambda^k \alpha_{32}(k\varphi)x_{02}+\hat\lambda^ka_{32}\xi^2_k+\\
&+\sum\limits_{j=0}^n\sum\limits_{i=0}^{j} \mu_{j,i} (y_{k1}-y_1^-)^{j-i} (y_{k2}-y_2^-)^{i}+\sum\limits_{i=0}^{n+1} A_i(y_{k1}-y_1^-)^{n+1-i}(y_{k2}-y_2^-)^{i}+c_3\hat\lambda^k\hat\xi_k+d_3\hat\gamma^{-k}\hat\eta_k+\ldots,\\
&\gamma^{-k}\cdot\left(-\overline y_{k1}\sin(k\psi)+\overline y_{k2}\cos(k\psi)\right)+\hat\gamma^{-k}\eta^2_k=\lambda^k \alpha_{41}(k\varphi)x_{01}+\hat\lambda^ka_{41}\xi^1_k+\lambda^k \alpha_{42}(k\varphi)x_{02}+\hat\lambda^ka_{42}\xi^2_k+\\
&+\sum\limits_{j=0}^n\sum\limits_{i=0}^{j} \nu_{j,i} (y_{k1}-y_1^-)^{j-i} (y_{k2}-y_2^-)^{i}+\sum\limits_{i=0}^{n+1}B_i (y_{k1}-y_1^-)^{n+1-i}(y_{k2}-y_2^-)^{i}+c_4\hat\lambda^k\hat\xi_k+d_4\hat\gamma^{-k}\hat\eta_k+\ldots,\\
&\overline u_0-u^+= \lambda^k \alpha_{51}(k\varphi)x_{01}+\hat\lambda^ka_{51}\xi^1_k+\lambda^k \alpha_{52}(k\varphi)x_{02}+\hat\lambda^ka_{52}\xi^2_k+ b_{51}(y_{k1}-y_1^-)+b_{52}(y_{k2}-y_2^-)+\\
&+c_5\hat\lambda^k\hat\xi_k+d_5\hat\gamma^{-k}\hat\eta_k+\ldots,\\
&v_k-v^-=\lambda^k \alpha_{61}(k\varphi)x_{01}+\hat\lambda^ka_{61}\xi^1_k+\lambda^k \alpha_{62}(k\varphi)x_{02}+\hat\lambda^ka_{62}\xi^2_k + b_{61}(y_{k1}-y_1^-)+b_{62}(y_{k2}-y_2^-)+\\
&+c_6\hat\lambda^k\hat\xi_k+d_6\hat\gamma^{-k}\hat\eta_k+\ldots,\\
\end{aligned}
\end{equation}
where for $i\in\{1,...,6\}$ the following equalities hold:
\begin{equation*}
\begin{aligned} 
&\alpha_{i1}(k\varphi)=a_{i1}\cos(k\varphi)+a_{i2}\sin(k\varphi), \;\;\;\;\;\;\;\;\;\;\;\;\;\;\;\;\;\;\; \alpha_{i2}(k\varphi)=a_{i2}\cos(k\varphi)-a_{i1}\sin(k\varphi).
\end{aligned}
\end{equation*}
Let us emphasize that in \eqref{lemmafirst-returnmapeq1} the dependence of the coefficients on the parameters is at least $C^{r-\max(2,n)}$-smooth since we lose two derivatives with respect to the parameters when write the local map in the form \eqref{eqlocal3} and $n$ derivatives when define $\mu_{j,i},\nu_{j,i}$.

Making the change of the variables:
\begin{equation*}
\begin{aligned} 
&X=x_0-x^+, \;\;\;\;\;\;\;\; Y=y_k-y^-, \;\;\;\;\;\;\;\; U=u_0-u^+,\;\;\;\;\;\;\;\; V= v_k-v^-,\\
\end{aligned}
\end{equation*}
and, taking into account the estimates on $\lambda,\hat\lambda,\gamma,\hat\gamma$ in Remark \ref{remark1}, we rewrite the system \eqref{lemmafirst-returnmapeq1} as
\begin{equation}\label{lemmafirst-returnmapeq2}
\begin{aligned} 
&\gamma^{-k}\cdot\left(\overline Y_1\cos(k\psi)+\overline Y_2\sin(k\psi)\right)=\sum\limits_{j=0}^n\sum\limits_{i=0}^{j} \hat\mu_{j,i} Y_1^{j-i} Y_2^{i}+\sum\limits_{i=0}^{n+1} A_i Y_1^{n+1-i}Y_2^{i}+h_1(X,Y,U,V,\overline X,\overline Y,\overline U,\overline V),\\
&\gamma^{-k}\cdot\left(-\overline Y_1\sin(k\psi)+\overline Y_2\cos(k\psi)\right)=\sum\limits_{j=0}^n\sum\limits_{i=0}^{j} \hat\nu_{j,i} Y_1^{j-i} Y_2^{i}+\sum\limits_{i=0}^{n+1} B_i Y_1^{n+1-i}Y_2^{i}+h_2(X,Y,U,V,\overline X,\overline Y,\overline U,\overline V),\\
&\left(\overline X_1,\overline X_2,\overline U,V\right)=h_3(X,Y,U,V,\overline X,\overline Y,\overline U,\overline V),
\end{aligned}
\end{equation}
where
\begin{equation*}
\begin{aligned} 
&\hat \mu_{0,0}=\mu_{0,0}-\gamma^{-k}\left(y_1^-\cos(k\psi)+y_2^-\sin(k\psi)\right)+O(\hat\gamma^{-k}), \;\;\;\;\;\;\;\;\;\;\;\;\;\; \hat\mu_{j,i}=\mu_{j,i}+O(\hat\gamma^{-k}),\\
&\hat \nu_{0,0}=\nu_{0,0}-\gamma^{-k}\left(-y_1^-\sin(k\psi)+y_2^-\cos(k\psi)\right)+O(\hat\gamma^{-k}), \;\;\;\;\;\;\;\;\;\;\;\; \hat\nu_{j,i}=\nu_{j,i}+O(\hat\gamma^{-k}),\\
\end{aligned}
\end{equation*}
and where the functions $h$ satisfy the estimates
\begin{equation*}
\begin{aligned} 
&h_{1,2}=o\left((Y_1^2+Y_2^2)^{\frac{n+1}{2}}\right)+O(\hat\gamma^{-k}) \;\;\;\; \text{and} \;\;\;\; h_3=O\left(|Y_1|+|Y_2|\right)+O(\hat\gamma^{-k}).\\
\end{aligned}
\end{equation*}

Expressing $\overline Y_1$, $\overline Y_2$ in the third and the fourth equations in \eqref{lemmafirst-returnmapeq2}, and applying the Implicit Function Theorem, we recast the system \eqref{lemmafirst-returnmapeq2} as
\begin{equation}\label{lemmafirst-returnmapeq3}
\begin{aligned} 
&\gamma^{-k}\overline Y_1=\sum\limits_{j=0}^n\sum\limits_{i=0}^{j} \tilde\mu_{j,i} Y_1^{j-i} Y_2^{i}+\sum\limits_{i=0}^{n+1}\tilde A_i Y_1^{n+1-i}Y_2^{i}+\tilde h_1(X,Y,U,\overline V),\\
&\gamma^{-k}\overline Y_2=\sum\limits_{j=0}^n\sum\limits_{i=0}^{j} \tilde\nu_{j,i} Y_1^{j-i} Y_2^{i}+\sum\limits_{i=0}^{n+1} \tilde B_i Y_1^{n+1-i}Y_2^{i}+ \tilde h_2(X,Y,U,\overline V),\\
&\left(\overline X_1,\overline X_2,\overline U,V\right)=\tilde h_3(X,Y,U,\overline V),\\
\end{aligned}
\end{equation}
where 
\begin{equation*}
\begin{aligned} 
&\tilde \mu_{j,i}=\hat\mu_{j,i}\cos(k\psi)-\hat\nu_{j,i}\sin(k\psi)+O(\hat\gamma^{-k}), \;\;\;\;\;\;\;\;\;\;\;\;\;\; \tilde \nu_{j,i}=\hat\mu_{j,i}\sin(k\psi)+\hat\nu_{j,i}\cos(k\psi)+O(\hat\gamma^{-k}),\\
&\tilde A_i=A_i\cos(k\psi)-B_i\sin(k\psi), \;\;\;\;\;\;\;\;\;\;\;\;\;\;\;\;\;\;\;\;\;\;\;\;\;\;\;\;\;\;\;\; \tilde B_i=A_i\sin(k\psi)+B_i\cos(k\psi),\\
\end{aligned}
\end{equation*}
and the functions $\tilde h$ satisfy 
\begin{equation*}
\begin{aligned} 
&\tilde h_{1,2}=o\left((Y_1^2+Y_2^2)^{\frac{n+1}{2}}\right)+O(\hat\gamma^{-k}) \;\;\;\; \text{and} \;\;\;\; \tilde h_3=O\left(|Y_1|+|Y_2|\right)+O(\gamma^k\cdot\hat\gamma^{-k}).\\
\end{aligned}
\end{equation*}

Further, let us nullify the constant terms in $\tilde h_3$ (by shifting the origin) and rescale the coordinates and the parameters as
\begin{equation}\label{scaling1}
\begin{aligned} 
&\left(Y_1,Y_2\right)=\gamma^{-\frac{k}{n}}\cdot\left(Y^{new}_1,Y^{new}_2\right),\;\;\;\;\;\;\;\;\;\; \left(X_1,X_2,U,V\right)=\frac{\gamma^{-\frac{k}{n}}}{\delta_k}\cdot\left(X^{new}_1,X^{new}_2,U^{new},V^{new}\right),\\
&\left(\tilde\mu_{j,i},\tilde\nu_{j,i}\right)=\gamma^{-k(1-\frac{j-1}{n})}\cdot\left(M_{j,i},N_{j,i}\right),\\
\end{aligned}
\end{equation}
where $\delta_{k}$ tends sufficiently slowly to zero as $k\rightarrow+\infty$. After it, the system \eqref{lemmafirst-returnmapeq3} is reduced to the desired form \eqref{lemmafirst-returnmapeq0}.

\end{demo}

\begin{remark}\label{remark2}
If instead if scaling \eqref{scaling1} we apply to the system \eqref{lemmafirst-returnmapeq3} the scaling
\begin{equation}\label{scaling2}
\begin{aligned} 
&\left(Y_1,Y_2\right)=\gamma^{-\frac{2k}{n}}\cdot\left(Y^{new}_1,Y^{new}_2\right),\;\;\;\;\;\;\;\;\;\; \left(X_1,X_2,U,V\right)=\frac{\gamma^{-\frac{2k}{n}}}{\delta_k}\cdot\left(X^{new}_1,X^{new}_2,U^{new},V^{new}\right),\\
&\left(\tilde\mu_{j,i},\tilde\nu_{j,i}\right)=\gamma^{-k(1-\frac{2(j-1)}{n})}\cdot\left(M_{j,i},N_{j,i}\right),\\
\end{aligned}
\end{equation}
then we obtain the following form for the first-return map:
\begin{equation}\label{remark2equation1}
\begin{aligned} 
&\overline Y_1=\sum\limits_{j=0}^n\sum\limits_{i=0}^{j} M_{j,i} Y_1^{j-i} Y_2^{i}+o(1)_{k\rightarrow+\infty},\\
&\overline Y_2=\sum\limits_{j=0}^n\sum\limits_{i=0}^{j} N_{j,i} Y_1^{j-i} Y_2^{i}+o(1)_{k\rightarrow+\infty},\\
&\left(\overline X_1,\overline X_2,\overline U,V\right)=o(1)_{k\rightarrow+\infty}, \\
\end{aligned}
\end{equation}
where the range of values of the new variables $(X,Y,U,V)$ and parameters $M, N$ covers a centered-at-zero ball whose radius tends to infinity as $k$ increases. The $o(1)$-terms tend to zero uniformly on any compact, along with the derivatives with respect to $(X,Y,U,V)$ up to the order $r$ and with respect to $M,N$ up to the order $\left(r-\max(2,n)\right)$.

\end{remark}

\begin{demo4}

Let us fix $k\in\mathbb N$, $k\ge 4$, and let $U^k$ be a unit ball in $\mathbb R^k$. Let us fix $r=2,\ldots,\infty,\omega$. Let us take any $\delta>0$ and any $C^r$ map $\zeta:U^k\rightarrow \mathbb R^k$ of the form 
\begin{equation}\label{approximationbypolynomial1}
\begin{aligned} 
\left(\overline X_1,\ldots,\overline X_{k-2}\right)&= 0, \\
\left(\overline X_{k-1},\overline X_k\right)&= \Phi \left(X_{k-1},X_k\right),\\
\end{aligned}
\end{equation}
where $\Phi: U^2\rightarrow \mathbb R^2$ is an arbitrary $C^r$ map. Let $\mathcal H(\zeta,\delta)$ be the set of all $k$-dimensional $C^r$ maps such that for each map from $\mathcal H(\zeta,\delta)$ there exists its renormalized iteration which is at distance smaller than $\delta$ from $\zeta$ (in the $C^r$ topology). By definition, the set $\mathcal H(\zeta,\delta)$ is open. Let us prove that this set is dense in the $ABR^*$-domain.

Let $f$ be any map belonging to the $ABR^*$-domain. By Theorem \ref{theorem2}, arbitrarily $C^r$-close to $f$ there exists a map $g$ which has infinitely many orbits of corank-2 tangency of every order. Let us consider one of such tangencies which is of order $n$, where $n=r$ if $r$ is finite and $n$ is sufficiently large if $r=\infty,\omega$. According to Remark \ref{remark2}, after $C^r$-small perturbation of this tangency, a first-return map in some coordinates is arbitrarily close to any map of the form
\begin{equation}\label{approximationbypolynomial2}
\begin{aligned} 
&\overline Y_1=\sum\limits_{j=0}^n\sum\limits_{i=0}^{j} M_{j,i} Y_1^{j-i} Y_2^{i},\;\;\;\;\;\;\; \overline Y_2=\sum\limits_{j=0}^n\sum\limits_{i=0}^{j} N_{j,i} Y_1^{j-i} Y_2^{i},\;\;\;\;\;\;\; \left(\overline X_1,\overline X_2,\overline U,V\right)=0. \\
\end{aligned}
\end{equation}
Therefore the Stone-Weierstrass Theorem in the smooth case and the Oka-Weil Theorem in the real-analytic case imply that for some $C^r$-small perturbation of the map $g$, its certain renormalized iteration approximates the map $\zeta$ with arbitrarily good accuracy. Thus, the set $\mathcal H(\zeta,\delta)$ is dense in the $ABR^*$-domain.

Let us take a sequence of maps $\left(\zeta_i\right)_{i\in\mathbb N}$ such that it is dense in the space of $C^r$ maps of the form \eqref{approximationbypolynomial1} and converging to zero sequence $\left(\delta_j\right)_{j\in\mathbb N}$, where $\delta_j$ is a positive real number. Then the intersection $\bigcap\limits_{i,j}\mathcal H(\zeta_i,\delta_j)$ is a residual set in the $ABR^*$-domain and consists of maps having $C^r$-universal two-dimensional dynamics. This completes the proof of the Theorem. 

\end{demo4}

\begin{center}
  \textbf{\Large{Acknowledgments}}
\end{center}

I am grateful to my scientific adviser Dmitry Turaev for setting the problem, his valuable comments and constant support. Also, I would like to thank M. Asaoka, P. Berger, S.V. Gonchenko, N. Gourmelon, V.P.H. Goverse, M. Helfter, A. Khodaeian Karim, J.S.W. Lamb, D. Li, S. van Strien for useful discussions. The research was supported by the EPSRC Centre for Doctoral Training in Mathematics of Random Systems: Analysis, Modelling and Simulation (EP/S023925/1), RSF grant (project 22-11-00027) and Leverhulme Trust grant RPG-2021-072.

\renewcommand\refname{References}

\end{document}